\documentclass[10pt,twocolumn,twoside,romanappendices]{IEEEtran}

\usepackage{amsmath,amssymb,amsthm}
\usepackage{graphicx}
\usepackage{algorithm,algpseudocode}
\usepackage{bm,color}
\usepackage{url}
\usepackage{cite}
\usepackage{textcomp}
\usepackage{stfloats}
\usepackage[breaklinks]{hyperref}
\usepackage[hyphenbreaks]{breakurl}

\allowdisplaybreaks[4]

\newtheorem{lemma}{Lemma}

\newtheorem{proposition}{Proposition}
\theoremstyle{definition}
\newtheorem{remark}{Remark}

\newcommand{\st}{\textrm{s.t.}}
\newcommand{\tV}{\textrm{v}}

\newcommand{\cA}{{\mathcal A}}
\newcommand{\cB}{{\mathcal B}}

\newcommand{\cP}{{\mathcal P}}

\newcommand{\cQ}{{\mathcal Q}}

\newcommand{\cM}{{\mathcal M}}
\newcommand{\cN}{{\mathcal N}}
\newcommand{\cJ}{{\mathcal J}}
\newcommand{\cW}{{\mathcal W}}
\newcommand{\cT}{{\mathcal T}}
\newcommand{\cV}{{\mathcal V}}

\DeclareMathOperator*{\argmax}{arg\,max}
\DeclareMathOperator*{\argmin}{arg\,min}

\begin{document}
\title{Robust Energy Management for Microgrids \\ With High-Penetration Renewables}%

\author{Yu Zhang,~\IEEEmembership{Student~Member,~IEEE,}~Nikolaos Gatsis,~\IEEEmembership{Member,~IEEE,}~and~Georgios B. Giannakis,~\IEEEmembership{Fellow,~IEEE}%
\thanks{%\protect\rule{0pt}{0.5cm}%
This work was supported by the University of Minnesota
Institute of Renewable Energy and the Environment (IREE) under Grant RL-0010-13.
This work was presented in part at the
\emph{3rd IEEE Intl. Conf. on Smart Grid Commun.}, Tainan, Taiwan, November 5--8, 2012.}
\thanks{
The authors are with the Dept. of ECE and the Digital Technology Center, University of Minnesota,
Minneapolis, MN 55455, USA. Tel/fax: +1(612)624-9510/625-2002. %
E-mails: {\tt \{zhan1220,gatsisn,georgios\}@umn.edu}  %
}
}

% The paper headers
\markboth{}
{}

\maketitle

% -------------------------------------------------------------------------
% Abstract
% -------------------------------------------------------------------------
\begin{abstract}
Due to its reduced communication overhead and robustness to
failures, distributed energy management is of paramount importance
in smart grids, especially in microgrids, which feature
distributed generation (DG) and distributed storage (DS).
Distributed economic dispatch for a microgrid with high renewable energy
penetration and demand-side management operating in
grid-connected mode is considered in this paper. To address the intrinsically stochastic availability of renewable
energy sources (RES), a novel power scheduling approach is introduced.
The approach involves the actual renewable energy as well as the energy
traded with the main grid, so that the supply-demand balance is maintained.
The optimal scheduling strategy minimizes the microgrid net cost,
which includes DG and DS costs, utility of dispatchable loads, and worst-case transaction cost
stemming from the uncertainty in RES. Leveraging the dual
decomposition, the optimization problem formulated is solved in a distributed fashion
by the local controllers of DG, DS, and dispatchable loads.
Numerical results are reported to corroborate the effectiveness of the novel approach.
\end{abstract}


\section*{Nomenclature}

\addcontentsline{toc}{section}{Nomenclature}

\subsection{Indices, numbers, and sets}

\begin{IEEEdescription}[\IEEEusemathlabelsep\IEEEsetlabelwidth{$M$, $m$}]

\item[$T$, $t$] Number of scheduling periods, period index.

\item[$M$, $m$] Number of conventional distributed generation (DG) units, and their index.

\item[$N$, $n$] Number of dispatchable (class-1) loads, load index.

\item[$Q$, $q$] Number of energy (class-2) loads, load index.

\item[$J$, $j$] Number of distributed storage (DS) units, and their index.

\item[$I$, $i$] Number of power production facilities with renewable energy source (RES), and facility index.

\item[$S$, $s$] Number of sub-horizons, and sub-horizon index.

\item[$k$] Algorithm iteration index.

\item[$\cT$] Set of time periods in the scheduling horizon.

\item[$\cT_s$] Sub-horizon $s$ for all RES facilities.

\item[$\cT_{i,s}$] Sub-horizon $s$ for RES facility $i$.

\item[$\cM$] Set of conventional DG units.

\item[$\cN$] Set of dispatchable loads.

\item[$\cQ$] Set of energy loads.

\item[$\cJ$] Set of DS units.

\item[$\cW$] Power output uncertainty set for all RES facilities.

\item[$\cW_{i}$] Power output uncertainty set of RES facility $i$.

\end{IEEEdescription}

\subsection{Constants}

\begin{IEEEdescription}[\IEEEusemathlabelsep\IEEEsetlabelwidth{$P_{D_n}^{\min,t}$, $P_{D_n}^{\max,t}$}]

\item[$P_{G_m}^{\min}$, $P_{G_m}^{\max}$] Minimum and maximum power output of conventional DG unit $m$.

%\item[$P_{G_m}^{\max}$] Maximum power output of conventional DG unit $m$.

\item[$R_{m,\text{up}}$, $R_{m,\text{down}}$] Ramp-up and ramp-down limits of conventional DG unit $m$.

%\item[$R_{m,\text{down}}$] Ramp-down rate of conventional DG unit $m$.

\item[$\mathsf{SR}^t$] Spinning reserve for conventional DG.

\item[$ L^t$] Fixed power demand of critical loads in period $t$.

\item[$P_{D_n}^{\min}$, $P_{D_n}^{\max}$] Minimum and maximum power consumption of load $n$.

%\item[$P_{D_n}^{\max,t}$] Maximum power consumption of  load $n$ in period $t$.

\item[$P_{E_q}^{\min,t}$, $P_{E_q}^{\max,t}$] Minimum and maximum power consumption of load $q$ in period $t$.

%\item[$P_{E_q}^{\max,t}$] Maximum power consumption of load $q$ in period $t$.

\item[$S_q$, $T_q$] Power consumption start and stop times of load $q$.

%\item[$T_q$] Power consumption end time of load $q$.

\item[$E_q^{\max}$] Total energy consumption of load $q$ from start time $S_q$ to termination time $T_q$.

\item[$P_{B_j}^{\min}$, $P_{B_j}^{\max}$] Minimum and maximum (dis)charging power of  DS unit $j$.

%\item[$P_{B_j}^{\max}$] Maximum (dis)charging power of storage $j$.

\item[$B_j^{\min}$] Minimum stored energy of DS unit $j$ in period $T$.

\item[$B_j^{\max}$] Capacity of DS unit $j$.

\item[$\eta_j$] Efficiency of DS unit $j$.

%\item[$\eta_j^{\textrm{ch}}$, $\eta_j^{\textrm{dis}}$] Charging and discharging efficiencies of DS unit $j$.

%\item[$ \eta_j^{\textrm{ch}}$] Charging efficiency of DS unit $j$.

\item[$P_{R}^{\min}$, $P_{R}^{\max}$] Lower and upper bounds for $P_{R}^t$.

%\item[$P_{R}^{\max}$] Upper limit for $P_{R}^t$.

\item[$\underline{W}_{i}^t$, $\overline{W}_i^t$] Minimum and maximum forecasted power output of RES facility $i$ in $t$.

%\item[$\overline{W}_i^t$] Maximum forecasted power output of RES facility $i$ in period $t$.

\item[$W^{\min}_{i,s}$, $W^{\max}_{i,s}$] Minimum and maximum forecasted total wind power of wind farm $i$ across sub-horizon $\cT_{i,s}$.

%\item[$W^{\max}_{i,s}$] Maximum forecasted total wind power of wind farm $i$ across time sub-horizon $\cT_{i,s}$.

\item[$W^{\min}_{s}$, $W^{\max}_{s}$] Minimum and maximum forecasted total wind power of all wind farms across sub-horizon $\cT_{s}$.

%\item[$W^{\max}_{s}$] Maximum forecasted total wind power of all the wind farms across time sub-horizon $\cT_{s}$.

\item[$\alpha^t$, $\beta^t$; $\gamma^t$, $\delta^t$] Purchase and selling prices; and functions thereof.

\item[$\pi_q^t$] Parameter of utility function of load $q$.

\item[$\mathsf{DOD}_j$; $\psi_j^t$] Depth of discharge specification of DS unit $j$; and parameters of storage cost.

\end{IEEEdescription}

\subsection{Uncertain quantities}

\begin{IEEEdescription}[\IEEEusemathlabelsep\IEEEsetlabelwidth{$P_{G_m}^t$}]

\item[$W_{i}^t$] Power output from RES facility $i$ in period $t$.

\end{IEEEdescription}

\subsection{Decision variables}

\begin{IEEEdescription}[\IEEEusemathlabelsep\IEEEsetlabelwidth{$\lambda^t$, $\mu^t$, $nu^t$}]

\item[$P_{G_m}^t$] Power output of DG unit $m$ in period $t$.

\item[$P_{D_n}^t$] Power consumption of load $n$ in period $t$.

\item[$P_{E_q}^t$] Power consumption of load $q$ in period $t$.

\item[$P_{B_j}^t$] (Dis)charging power of DS unit $j$ in period $t$.

\item[$B_j^t$] Stored energy of DS unit $j$ at the end of the period $t$.

\item[$P_R^t$] Net power delivered to the microgrid from the RES and storage in period $t$.

\item[$\tilde{P}_R^t$] Auxiliary variable.

\item[$\mathbf{x}$] Vector collecting all decision variables.

\item[$\lambda^t$, $\mu^t$, $\nu^t$] Lagrange multipliers.

\item[$\mathbf{z}$] Vector collecting all Lagrange multipliers.

\item[$W_{\textrm{worst}}^t$] Power production from all RES facilities in $t$ yielding the worst-case transaction cost.

\end{IEEEdescription}

\subsection{Functions}

\begin{IEEEdescription}[\IEEEusemathlabelsep\IEEEsetlabelwidth{$G(\cdot)$, $\tilde G(\cdot)$}]

\item[$C_m^t(\cdot)$] Cost of conventional DG unit $m$ in period $t$.

\item[$U_{D_n}^t(\cdot)$] Utility of load $n$ in period $t$.

\item[$U_{E_q}^t(\cdot)$] Utility of load $q$ in period $t$.

\item[$H_j^t(\cdot)$] Cost of DS unit $j$ in period $t$.

\item[$G(\cdot,\cdot)$] Worst transaction cost across entire horizon.

\item[$G(\cdot)$, $\tilde G(\cdot)$] Modified worst-case transaction cost.

\item[$\mathcal{L}(\mathbf{x},\mathbf{z})$] Lagrangian function.

\item[$\mathcal{D}(\mathbf{z})$] Dual function.

\end{IEEEdescription}

\section{Introduction}
Microgrids are power systems comprising distributed energy resources
(DERs) and electricity end-users, possibly with controllable elastic
loads, all deployed across a limited geographic
area~\cite{Hatziargyriou-PESMag}. Depending on their origin, DERs
can come either from distributed generation (DG) or from distributed
storage (DS). DG refers to small-scale power generators such as diesel generators,
fuel cells, and renewable energy sources (RES), as in wind or photovoltaic (PV)
generation. DS paradigms include batteries, flywheels, and pumped
storage. Specifically, DG brings power closer to the point it is
consumed, thereby incurring fewer thermal losses and bypassing
limitations imposed by a congested transmission network. Moreover,
the increasing tendency towards high penetration of RES stems from
their environment-friendly and price-competitive advantages over
conventional generation. Typical microgrid loads include critical
non-dispatchable types and elastic controllable ones.

\begin{figure}
\centering
\includegraphics[scale=0.55]{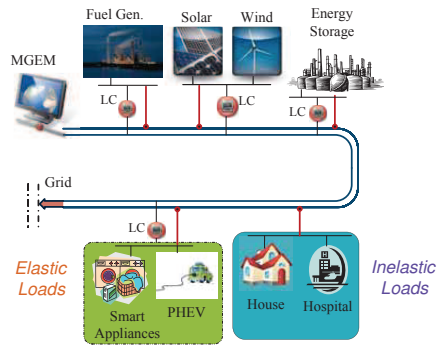}
\caption{Distributed control and computation architecture of a
microgrid.}
%system. The microgrid energy manager (MGEM) coordinates
%the local controllers (LCs) of DERs and dispatchable loads. }
\label{fig:MGModel}
\vspace*{-0.5cm}
\end{figure}

Microgrids operate in grid-connected or island mode, and may entail distribution networks with residential or
commercial end-users, in rural or urban areas. A typical
configuration is depicted in Fig.~\ref{fig:MGModel}; see
also~\cite{Hatziargyriou-PESMag}.
The microgrid energy manager (MGEM) coordinates the DERs and the controllable loads.
Each of the DERs and loads has a local controller (LC), which coordinates with
the MGEM the scheduling of resources through the communications
infrastructure in a distributed fashion.
The main challenge in energy scheduling is to account for the random and nondispatchable nature of the RES.

Optimal energy management for microgrids including
economic dispatch (ED), unit commitment (UC), and demand-side management (DSM)
is addressed in~\cite{Stluka11}, but without pursuing a robust formulation against RES uncertainty.
Based on the Weibull distribution for wind speed and the wind-speed-to-power-output mappings, an ED problem
is formulated to minimize the risk of overestimation and underestimation of available wind power~\cite{HetzerYB08}.
%Lyapunov stochastic optimization has been adopted to maximize the long-term profit of a RES facility in~\cite{Neely10}.
Stochastic programming is also used to cope with the variability of RES. Single-period
chance-constrained ED problems for RES have been studied
in~\cite{LiuX10}, yielding probabilistic guarantees that the load
will be served. Considering the uncertainties of demand profiles
and PV generation, a stochastic program is formulated to minimize the overall
cost of electricity and natural gas for a building in~\cite{Guan10}.
Without DSM, robust scheduling problems with penalty-based costs for uncertain supply
and demand have been investigated in~\cite{BertsimasLSZZ10}. Recent works explore energy
scheduling with DSM and RES using only centralized algorithms~\cite{JiangL11-Allerton,ZhaoZ10}.
An energy source control and DS planning problem for a microgrid
is formulated and solved using model predictive control in~\cite{Jin11}.
%Using energy storage to mitigate fluctuation in generation due to time-varying RES,
%an optimal power flow problem is formulated in~\cite{Chandy10}.
Distributed algorithms are developed in~\cite{DGH-CDC}, but they only
coordinate DERs to supply a given load without considering the stochastic nature of RES.
%In all the aforementioned works however, robust formulations accounting
%for the RES randomness are not pursued.
Recently, a worst-case transaction cost based
energy scheduling scheme has been proposed to address the variability of RESs through
robust optimization that can also afford distributed implementation~\cite{YuSGComm12}.
However,~\cite{YuSGComm12} considers only a single wind farm and no DS, and
its approach cannot be readily extended to include multiple RESs and DS.

The present paper deals with optimal energy
management for both supply and demand of a grid-connected microgrid incorporating
RES. The objective of minimizing the microgrid net cost accounts for
conventional DG cost, utility of elastic loads, penalized cost of DS, and a worst-case
transaction cost. The latter stems from the ability of the microgrid
to sell excess energy to the main grid, or to import energy in case
of shortage. A \emph{robust} formulation accounting for the
worst-case amount of harvested RES is developed. A novel model is
introduced in order to maintain the supply-demand balance arising
from the intermittent RES. Moreover, a transaction-price-based
condition is established to ensure convexity of the overall problem (Section~\ref{sec:Probformulation}).
The separable structure and strong duality of the resultant
problem are leveraged to develop a low-overhead
\emph{distributed} algorithm based on dual decomposition, which is
computationally efficient and resilient to communication outages or attacks.
%The distributed implementation relies upon message exchanges between the MGEM and LCs.
For faster convergence, the proximal bundle method is employed for the non-smooth
subproblem handled by the LC of RES (Section~\ref{Section:Dualmethod}).
Numerical results corroborate the merits of the novel designs (Section~\ref{sec:Numericalresults}),
and the paper is wrapped up with a concluding summary (Section~\ref{sec:Conclusions}).

Compared to~\cite{YuSGComm12}, the contribution of the paper is threefold, and
of critical importance for microgrids with high-penetration renewables.
First, a detailed model for DS is included, and different design choices for storage cost functions are given
to accommodate, for example, depth-of-discharge specifications.
Second, with the envisioned tide of high-penetration renewable energy, multiple
wind farms are considered alongside two pertinent uncertainty models.
Finally, a new class of controllable loads is added, with each load having a requirement
of total energy over the scheduling horizon, as is the case with charging of plug-in hybrid electric vehicles (PHEVs).
Detailed numerical tests are presented to illustrate the merits of the scheduling decisions for the DG, DS, and controllable loads.

\noindent \emph{Notation}. Boldface lower case letters represent
vectors; $\mathbb{R}^{n}$ and $\mathbb{R}$ stand for spaces of $n
\times 1$ vectors and real numbers, respectively;
$\mathbb{R}^{n}_{+}$ is the $n$-dimensional non-negative orthant;
$\mathbf{x}^{\prime}$ transpose, and $\|\mathbf{x}\|$ the Euclidean
norm of $\mathbf{x}$.

\section{Robust Energy Management Formulation}\label{sec:Probformulation}
%Consider a microgrid comprising $M$ conventional (fossil fuel) generators, $N$ controllable (dispatchable) loads, $I$ RES facilities, and $J$ distributed storage devices (see also Fig.~\ref{fig:MGModel}).
%The scheduling horizon is $\cT:=\{1,2,\ldots,T\}$ (e.g., one-day ahead).
%Let $P_{G_m}^t$ be the power produced by the $m$th conventional generator, and $P_{D_n}^t$ the power consumed by the $n$th dispatchable load,
%where $m \in \cM := \{1,\ldots,M\}$, $n \in \cN :=\{1,\ldots,N\}$, and $t \in \cT$.
%The power delivered to the microgrid to keep the supply-demand balance at slot $t$ is denoted by $P_R^t$.
%The ensuing subsection details the storage model, the RES uncertainty model, and the transaction mechanism between the microgrid and the main grid. Subsection~\ref{Subsec:Form} formulates the microgrid energy management problem, which boils down to optimally scheduling the variables $P_{G_m}^t$, $P_{D_n}^t$, and $P_R^t$.

Consider a microgrid comprising $M$ conventional (fossil fuel)
generators, $I$ RES facilities, and $J$ DS units (see also
Fig.~\ref{fig:MGModel}). The scheduling horizon is
$\cT:=\{1,2,\ldots,T\}$ (e.g., one-day ahead). The particulars of
the optimal scheduling problem are explained in the next
subsections.

%The ensuing subsection details the storage model, the RES uncertainty model, and the transaction mechanism between the microgrid and the main grid. Subsection~\ref{Subsec:Form} formulates the microgrid energy management problem, which boils down to optimally scheduling the variables.

\subsection{Load Demand Model}
%For the microgrid power end-users, the loads may have diverse usage features. Generally, they can be
Loads are classified in two categories. The first comprises
inelastic loads, whose power demand should be satisfied at all
times. Examples are power requirements of hospitals or illumination
demand from residential areas.
%The total fixed power demand from critical loads at slot $t$ is denoted by $L^t$.

The second category consists of elastic loads, which are
dispatchable, in the sense that their power consumption is
adjustable, and can be scheduled.
%which are dispatchable according to the real-time price to reduce the peak demand for electricity.
These loads can be further divided in two classes, each having the following
characteristics:
\begin{itemize}
\item[i)] The first class contains loads with power consumption
$P_{D_n}^t \in [P_{D_n}^{\min}, P_{D_n}^{\max}]$, where $n \in \cN
:=\{1,\ldots,N\}$, and $t \in \cT$. Higher power consumption yields
higher utility for the end user. The utility function of the $n$th
dispatchable load, $U_{D_n}^t(P_{D_n}^t)$, is selected to be increasing and concave,
with typical choices being piecewise linear or smooth quadratic; see also~\cite{ChenDoyle-2market}.
An example from this class is an A/C.
\item[ii)] The second class includes loads indexed by $q \in \cQ:=\{1,\ldots,Q\}$
with power consumption limits $P_{E_q}^{\min}$ and $P_{E_q}^{\max}$,
and prescribed total energy requirements $E_q$ which have to be
achieved from the start time $S_q$ to termination time $T_q$; see e.g.,~\cite{RadWJSL10}.
This type of loads can be the plug-in hybrid electric vehicles (PHEVs).
Power demand variables $\{P_{E_q}^t\}_{t=1}^{T}$ therefore are
constrained as $\sum_{t=S_q}^{T_q}P_{E_q}^{t} = E_q$ and
$P_{E_q}^t \in [P_{E_q}^{\min,t}, P_{E_q}^{\max,t}],~t \in \cT$, while
$P_{E_q}^{\min,t}=P_{E_q}^{\max,t}=0$ for $t \notin \{S_q,\ldots,T_q\}$.
%a slab $\cP_{E_q} := \{\{P_{E_q}^t\}_{t=1}^{T}|\sum_{t \in
%\cS_q}P_{E_q}^{t} = E_q, P_{E_q}^{\min,t} \le P_{E_q}^t \le
%P_{E_q}^{\max,t},~t \in \cT\}$, where $\cS_q :=\{S_q,\ldots,T_q\}$
%and $P_{E_q}^{\min,t}= P_{E_q}^{\max,t}=0$ for $t
%\notin \cS_q$.
Higher power consumption in earlier slots as opposed
to later slots may be desirable for a certain load, so that the
associated task finishes earlier. This behavior can be encouraged by
adopting for the $q$th load an appropriately designed time-varying concave utility function
$U_{E_q}^t(P_{E_q}^t)$. An example is $U_{E_q}^t(P_{E_q}^t) :=
\pi_q^tP_{E_q}^t$, with  weights $\{\pi_q^t\}$ decreasing in $t$
from slots $S_q$ to $T_q$. Naturally, $U_{E_q}^t(P_{E_q}^t) \equiv 0$ can be selected
if the consumer is indifferent to how power is consumed across slots.
%Therefore, this weighted sum utility encourages scheduling more power
%to the initially several slots of $\tau_q$, e.g., meaning try to charge the PHEV as soon as possible.
%\textcolor{red}{(other two examples for utility functions, one in proposal, another is the length vector)}
\end{itemize}

\subsection{Distributed Storage Model}
Let $B_j^t$ denote the stored energy of the $j$th battery
at the end of the slot $t$,
with initial available energy $B_j^0$ while $B_j^{\max}$ denotes the
battery capacity, so that $0 \le B_j^t \le
B_j^{\max},~j \in \cJ := \{1,\ldots,J\}$. Let $P_{B_j}^t$ be the power delivered to
(drawn from) the $j$th storage device at slot $t$, which amounts to
charging ($P_{B_j}^t\geq 0$) or discharging ($P_{B_j}^t\leq 0$) of
the battery. Clearly, the stored energy obeys the dynamic equation
\begin{align}
B_j^{t} = B_j^{t-1} + P_{B_j}^t,~j \in \cJ,~t \in \cT \;.
\end{align}
Variables $P_{B_j}^t$ are constrained in the following ways:

\begin{itemize}
\item[i)] The amount of (dis)charging is bounded, that is
\begin{align}
P_{B_j}^{\min} \le &P_{B_j}^t \le P_{B_j}^{\max} \\
-\eta_j B_j^{t-1} \le &P_{B_j}^t  % P_{B_j}^t %\le (B_j^{\max}-B_j^{t-1})/\eta_j^{\textrm{ch}}
\label{EffConstr}
\end{align}
with bounds $P_{B_j}^{\min}<0$ and $P_{B_j}^{\max}>0$, while
$\eta_j\in(0,1]$ is the efficiency of DS unit $j$ \cite{aamas,Alimisis12}.
The constraint in~\eqref{EffConstr} means that a fraction $\eta_j$ of the stored energy $B_j^{t-1}$ is available for discharge.
%The second constraint means that the energy required to fully charge the battery at slot $t$ is $(B_j^{\mathrm{max}}-B_j^{t-1})/\eta_j^{\mathrm{ch}}$.}

\item[ii)] Final stored energy is also bounded for the sake of future scheduling horizons,
that is $B_j^T \ge B_j^{\min}$.
\end{itemize}

%The storage cost $H_j^t(B_j^t)$ can be also taken into account.
%Considering the battery reserve, a penalty can be postulated
%proportional to the deviation from the battery capacity as
%$H_j^t(B_j^t) := \varpi_j^t(B_j^{\textrm{max}}-B_j^t)$, meaning it
%is desirable to maintain a close-to-full battery level at the end of
%each period\cite{Chandy10}.
To maximize DS lifetime, a storage cost $H_j^t(B_j^t)$ can be employed to encourage the stored energy
to remain above a specified depth of discharge, denoted as $\mathsf{DOD}_j\in [0,1]$,
where $100\%$ $(0\%)$ depth of discharge means the battery is empty (full)~\cite{Alimisis12}.
Such a cost is defined
as $H_j^t(B_j^t) := \psi_j^t[(1-\mathsf{DOD}_j)B_j^{\textrm{max}}-B_j^t]$.
Note that the storage cost $H_j^t(B_j^t)$ can be interpreted as imposing a soft constraint
preventing large variations of the stored energy.
Clearly, %weights $\{\psi_j^t\}$ control the variation of the stored energy;
higher weights $\{\psi_j^t\}$ encourage smaller variation.
If high power exchange is to be allowed,
these weights can be chosen very small, or one can even select $H_j^t(B_j^t)\equiv 0$ altogether.

\subsection{Worst-case Transaction Cost}
Let $W_i^t$ denote the \emph{actual} renewable energy harvested by
the $i$th RES facility at time slot $t$, and also let $\mathbf{w}$ collect all $W_i^t$, i.e.,
$\mathbf{w} := [W_1^1,\ldots,W_1^T,\ldots,W_I^1,\ldots,W_I^T]$.
To capture the intrinsically stochastic and time-varying availability of RES,
it is postulated that $\mathbf{w}$ is unknown, but lies in a polyhedral uncertainty set $\cW$.
The following are two practical examples.
\begin{itemize}
\item[i)] The first example postulates a separate uncertainty set $\cW_i$ for each RES facility in the form
\begin{align}
\label{perW}
\cW_i &:=
\Bigg\{\{W_i^t\}_{t=1}^{T}|\underline{W}_{i}^t \le W_i^t \le \overline{W}_i^t, \nonumber \\
&W^{\min}_{i,s} \le \sum\limits_{t \in \cT_{i,s}} W_i^t\le W^{\max}_{i,s},
\cT = \bigcup \limits _{s=1}^{S} \cT_{i,s} \Bigg\}
\end{align}
where $\underline{W}_{i}^t$ ($\overline{W}_i^t$) denotes a lower
(upper) bound on $W_i^t$; $\cT$ is partitioned
into consecutive but non-overlapping sub-horizons $\cT_{i,s}$ for $i = 1,\ldots,I$,
$s=1,2,\ldots,S$; the total renewable energy for the $i$th RES facility over the $s$th
sub-horizon is assumed bounded by $W^{\min}_{i,s}$ and $W^{\max}_{i,s}$.
In this example, $\cW$ takes the form of Cartesian product
\begin{align}
\cW = \cW_1 \times\ldots\times\cW_I.
\end{align}
\item[ii)] The second example assumes a joint uncertainty model across all the RES facilities as
\begin{align}
\label{jointW}
\cW &:=
\Bigg\{\mathbf{w}| \underline{W}_{i}^t \le W_i^t \le \overline{W}_i^t, \nonumber \\
&W^{\min}_{s} \le \sum\limits_{t \in \cT_s}\sum\limits_{i=1}^{I}
W_i^t\le W^{\max}_{s}, \cT = \bigcup \limits _{s=1}^{S} \cT_s\Bigg\}
\end{align}
where $\underline{W}_{i}^t$ ($\overline{W}_i^t$) denotes a lower
(upper) bound on $W_i^t$; $\cT$ is partitioned
into consecutive but non-overlapping sub-horizons $\cT_s$ for $s=1,2,\ldots,S$;
the total renewable energy harvested by all the RES facilities over the $s$th sub-horizon is bounded by
$W^{\min}_{s}$ and $W^{\max}_{s}$; see also~\cite{ZhaoZ10}.
\end{itemize}

The previous two RES uncertainty models are quite general and
can take into account different geographical and meteorological
factors. The only information required is the deterministic lower
and upper bounds, namely $\underline{W}_{i}^t$, $\overline{W}_i^t$,
$W^{\min}_{i,s}$, $W^{\max}_{i,s}$, $W^{\min}_{s}$, $W^{\max}_{s}$,
which can be determined via inference schemes based
on historical data~\cite{Pinson10}.

Supposing the microgrid operates in a grid-connected mode, a
transaction mechanism between the microgrid and the main grid is
present, whereby the microgrid can buy/sell energy from/to the spot
market. Let $P_R^t$ be an auxiliary variable denoting the net power
delivered to the microgrid from the renewable energy sources and the
distributed storage in order to maintain the supply-demand balance
at slot $t$. The shortage energy per slot $t$ is given by
$\left[P_R^t-\sum_{i=1}^{I}W_i^t + \sum_{j=1}^{J}P_{B_j}^t\right]^{+}$,
while the surplus energy is
$\left[P_R^t-\sum_{i=1}^{I}W_i^t + \sum_{j=1}^{J}P_{B_j}^t\right]^{-}$,
where $[a]^{+}:=\max\{a,0\}$, and $[a]^{-}:=\max\{-a,0\}$.

The amount of shortage energy is bought with known purchase price
$\alpha^t$, while the surplus energy is sold to the main grid with
known selling price $\beta^t$. The worst-case net transaction cost
is thus given by
\begin{align}
G(\{P_R^t\},\{P_{B_j}^t\}) := \mathop{\max}\limits_{\mathbf{w} \in \cW}
&\sum_{t=1}^T \Bigg(\alpha^t\Bigg[P_R^t-\sum_{i=1}^{I}W_i^t + \sum_{j=1}^{J}P_{B_j}^t\Bigg]^{+} \nonumber \\
&\hspace{-1.5cm}-\beta^t\Bigg[P_R^t-\sum_{i=1}^{I}W_i^t + \sum_{j=1}^{J}P_{B_j}^t\Bigg]^{-}\Bigg)
\end{align}
where $\{P_R^t\}$ collects $P_R^t$ for $t = 1,2,\ldots,T$
and $\{P_{B_j}^t\}$ collects $P_{B_j}^t$ for $j = 1,2,\ldots,J,~t = 1,2,\ldots,T$.
\begin{remark}\textit{(Worst-case model versus stochastic model)}.
The worst-case robust model advocated here is particularly attractive when the
probability distribution of the renewable power production is unavailable.
This is for instance the case for multiple wind farms,
where the spatio-temporal joint distribution
of the wind power generation is intractable (see detailed discussions in~\cite{Yu-isgt13} and~\cite{Morales10}).
If an accurate probabilistic model is available, an expectation-based stochastic program can be formulated to bypass
the conservatism of worst-case optimization.
In the case of wind generation, suppose that wind power $W_i^t$ is a function of the random wind velocity $v_i^t$,
for which different models are available, and the wind-speed-to-power-output mappings $W_i^t(v_i^t)$ are known~\cite{CartaRV09}.
Then, the worst-case transaction cost can be replaced by the \emph{expected} transaction cost
$G(\{P_R^t\},\{P_{B_j}^t\}):=\mathbb{E}_{\mathbf{v}}\Big(\sum_{t=1}^T \alpha^t
[P_R^t-\sum_{i=1}^{I}W_i^t(v_i^t)+\sum_{j=1}^{J}P_{B_j}^t]^{+} -
\beta^t[P_R^t-\sum_{i=1}^{I}W_i^t(v_i^t)+\sum_{j=1}^{J}P_{B_j}^t]^{-}\Big)$,
where $\mathbf{v}$ collects $v_i^t$ for all $i$ and $t$.
\end{remark}

\subsection{Microgrid Energy Management Problem}
\label{Subsec:Form}

Apart from RES, microgrids typically entail also conventional DG.
Let $P_{G_m}^t$ be the power produced by the $m$th conventional
generator, where $m \in \cM := \{1,\ldots,M\}$ and $t \in \cT$. The
cost of the $m$th generator is given by an increasing convex function $C_m^t(P_{G_m}^t)$,
which typically is either piecewise linear or smooth quadratic.

The energy management problem amounts to minimizing the microgrid
social net cost; that is, the cost of conventional generation,
storage, and the worst-case transaction cost (due to the
volatility of RES) minus the utility of dispatchable loads:
\begin{subequations}
\label{DERsched}
\begin{align}
\textrm{(P1)}
&\mathop{\min}_{\substack{\{P_{G_m}^t,P_{D_n}^t, \\ P_{E_q}^t,B_j^t,P_{B_j}^t,P_R^t\}}}
\sum^{T}_{t=1}\Bigg(\sum^{M}_{m=1}C_m^t(P_{G_m}^t)-\sum^{N}_{n=1}U_{D_n}^t(P_{D_n}^t) \nonumber \\
&\hspace{-0.7cm} -\sum^{Q}_{q=1}U_{E_q}^t(P_{E_q}^t)+ \sum^{J}_{j=1}H_j^t(B_j^t)\Bigg)+G(\{P_R^t\},\{P_{B_j}^t\}) \label{ObjFunc}\\
& \textrm{subject to:}  \nonumber \\
&P_{G_m}^{\min} \le P_{G_m}^t \le P_{G_m}^{\max},~m \in \cM, ~t \in \cT\label{Plimits}\\
&P_{G_m}^t-P_{G_m}^{t-1} \le R_{m,\text{up}},~m \in \cM, ~t \in \cT\label{RampUp}\\
&P_{G_m}^{t-1}-P_{G_m}^t \le R_{m,\text{down}},~m \in \cM, ~t \in \cT\label{RampDown}\\
&\sum^{M}_{m=1}(P_{G_m}^{\max}-P_{G_m}^{t}) \ge \mathsf{SR}^t,~t \in \cT\label{SpinRes}\\
&P_{D_n}^{\min} \le P_{D_n}^t \le P_{D_n}^{\max},~n \in \cN, ~t \in \cT\label{Dlimits}\\
&P_{E_q}^{\min,t} \le P_{E_q}^t \le P_{E_q}^{\max,t},~q \in \cQ, ~t \in \cT\label{Elimits}\\
&\sum^{T_q}_{t=S_q}P_{E_q}^{t} = E_q,~q \in \cQ \label{Esumlimits}\\
&0 \le B_j^t \le B_j^{\max},~B_j^T \ge B_j^{\min},~j \in \cJ, ~t \in \cT\label{Blimits}\\
&P_{B_j}^{\min} \le P_{B_j}^t \le P_{B_j}^{\max},~j \in \cJ, ~t \in \cT\label{PBlimits1}\\
& -\eta_j B_j^{t-1} \le P_{B_j}^t,~j \in \cJ, ~t \in \cT\label{PBlimits2}\\
&B_j^{t} = B_j^{t-1}+ P_{B_j}^t,~j \in \cJ, ~t \in \cT\label{ChargingEq}\\
%&B_j^T \ge B_j^{\min},~j \in \cJ \label{FinalSoC}\\
&P_{R}^{\min} \le P_R^t \le P_{R}^{\max},~t \in \cT\label{Rlimits}\\
&\hspace{-0.7cm}  \sum^{M}_{m=1}P_{G_m}^t+P_R^t = L^t+\sum^{N}_{n=1}P_{D_n}^t+\sum^{Q}_{q=1}P_{E_q}^t,~t \in \cT.\label{Balance}
\end{align}
\end{subequations}

Constraints \eqref{Plimits}--\eqref{SpinRes} stand for the
minimum/maximum power output, ramping up/down limits, and spinning
reserves, respectively, which capture the typical physical
requirements of a power generation system. Constraints
\eqref{Dlimits} and \eqref{Rlimits} correspond to the
minimum/maximum power of the flexible load demand and committed
renewable energy. Constraint~\eqref{Balance} is the power
supply-demand \emph{balance equation} ensuring the total demand is
satisfied by the power generation at any time.

Note that constraints~\eqref{Plimits}--\eqref{Balance} are linear,
while $C_m^t(\cdot)$, $-U_{D_n}^t(\cdot)$, $-U_{E_q}^t(\cdot)$, and $H_j^t(\cdot)$
are convex (possibly non-differentiable or
non-strictly convex) functions. Consequently, the convexity of (P1)
depends on that of $G(\{P_R^t\},\{P_{B_j}^t\})$, which is established in
the following proposition.
\begin{proposition}
\label{prop:convex} If the selling price $\beta^t$ does not exceed
the purchase price $\alpha^t$ for any $t \in \mathcal{T}$, then the
worst-case transaction cost $G(\{P_R^t\},\{P_{B_j}^t\})$ is convex in
$\{P_R^t\}$ and $\{P_{B_j}^t\}$.
\end{proposition}
\begin{IEEEproof}
Using that $[a]^{+}+[a]^{-}=|a|$, and $[a]^{+}-[a]^{-}=a$,
$G(\{P_R^t\},\{P_{B_j}^t\})$ can be re-written as
\begin{align}
G(\{P_R^t\},\{P_{B_j}^t\})=\mathop{\max}\limits_{\mathbf{w} \in \cW}
&\sum_{t=1}^T\Bigg(\delta^t\Bigg|P_R^t-\sum_{i=1}^{I}W_i^t + \sum_{j=1}^{J}P_{B_j}^t\Bigg| \nonumber \\
&\hspace{-1.5cm}+\gamma^t\Bigg(P_R^t-\sum_{i=1}^{I}W_i^t + \sum_{j=1}^{J}P_{B_j}^t\Bigg)\Bigg)
\end{align}
with $\delta^t := (\alpha^t-\beta^t)/2$, and $\gamma^t :=(\alpha^t+\beta^t)/2$.
Since the absolute value function is convex,
and the operations of nonnegative weighted summation and pointwise
maximum (over an infinite set) preserve
convexity~\cite[Sec.~3.2]{Boyd}, the claim follows readily.
\end{IEEEproof}

An immediate corollary of Proposition~\ref{prop:convex} is that the
energy management problem~(P1) is convex if $\beta^t\leq \alpha^t$
for all $t$. The next section focuses on this case, and designs an
efficient decentralized solver for (P1).

\section{Distributed Algorithm}\label{Section:Dualmethod}
In order to facilitate a distributed algorithm for (P1), a variable
transformation is useful. Specifically, upon introducing
$\tilde{P}_R^t:=P_R^t+\sum_{j=1}^JP_{B_j}^t$, (P1) can be re-written
as
\begin{subequations}
\begin{align}
\hspace{-0.25cm} \textrm{(P2)} \quad
&\mathop{\min}_{\mathbf{x}}
\sum^{T}_{t=1}\Bigg(\sum^{M}_{m=1}C_m^t(P_{G_m}^t)-\sum^{N}_{n=1}U_{D_n}^t(P_{D_n}^t) \nonumber \\
&-\sum^{Q}_{q=1}U_{E_q}^t(P_{E_q}^t)+ \sum^{J}_{j=1}H_j^t(B_j^t)\Bigg)+G(\{\tilde{P}_R^t\}) \label{P2ObjFunc}\\
&\textrm{subject to:} \quad \eqref{Plimits}-\eqref{Balance}\nonumber \\
&\tilde{P}_R^t=P_R^t+\sum_{j=1}^JP_{B_j}^t,~t \in \cT\label{VarChange}
\end{align}
\end{subequations}
where $\mathbf{x}$ collects all the primal variables $\{P_{G_m}^t,P_{D_n}^t,P_{E_q}^t,P_{B_j}^t,B_j^t,P_R^t,\tilde{P}_R^t\}$;
$\{\tilde{P}_R^t\}$ collects $\tilde{P}_R^t$ for $t=1,\ldots,T$; and {\small
\begin{align}
G(\{\tilde{P}_R^t\}) :=\mathop{\max}\limits_{\mathbf{w} \in \cW}
\sum_{t=1}^T\left(\delta^t\left|\tilde{P}_R^t-\sum_{i=1}^{I}W_i^t\right|
+\gamma^t\left(\tilde{P}_R^t-\sum_{i=1}^{I}W_i^t\right)\right).
\end{align}
}
The following proposition extends the result of Proposition~\ref{prop:convex}
to the transformed problem, and asserts its strong duality.

%% long formulas span two columns
\newcounter{TempEqCnt}
\setcounter{TempEqCnt}{\value{equation}}
\setcounter{equation}{15}
\begin{figure*}[hb]
\centering
\hrulefill
\begin{align}
&\{P_{G_m}^{t}(k)\}_{t=1}^T \in \argmin\limits_{\substack{\{P_{G_m}^t\} \\ \st~\eqref{Plimits}-\eqref{RampDown}}}
\Bigg\{\sum^{T}_{t=1}\bigg(C_m^t(P_{G_m}^t)+\big(\mu^t(k)-\lambda^t(k)\big)P_{G_m}^t\bigg)\Bigg\} \label{PGsubp}\\
&\{P_{D_n}^t(k)\}_{t=1}^T  \in \argmin\limits_{\substack{\{P_{D_n}^t\} \\ \st~\eqref{Dlimits}}}
\Bigg\{\sum^{T}_{t=1}\bigg(\lambda^t(k)P_{D_n}^t-U_{D_n}^t(P_{D_n}^t)\bigg)\Bigg\} \label{PDsubp} \\
&\{P_{E_q}^t(k)\}_{t=1}^T \in \argmin\limits_{\substack{\{P_{E_q}^t\} \\ \st~\eqref{Elimits}-\eqref{Esumlimits}}}
\Bigg\{\sum^{T}_{t=1}\bigg(\lambda^t(k)P_{E_q}^t-U_{E_q}^t(P_{E_q}^t)\bigg)\Bigg\} \label{PEsubp} \\
&\{P_{B_j}^t(k)\}_{t=1}^T \in \argmin\limits_{\substack{\{P_{B_j}^t,B_j^t\} \\ \st~\eqref{Blimits}-\eqref{ChargingEq}}}
\Bigg\{\sum^{T}_{t=1}\bigg(\nu^t(k)P_{B_j}^t+H_j^t(B_j^t)\bigg)\Bigg\} \label{PBsubp} \\
&\{P_R^t(k),\tilde{P}_R^t(k)\}_{t=1}^T \in \argmin\limits_{\substack{\{P_{R}^t,\tilde{P}_R^t\} \\ \st~\eqref{Rlimits}}}
\Bigg\{\sum^{T}_{t=1}\bigg(\big(\nu^t(k)-\lambda^t(k)\big)P_R^t\bigg)
 +G(\{\tilde{P}_R^t\})-\sum^{T}_{t=1}\nu^t(k)\tilde{P}_R^t\Bigg\}\label{PRsubp}
\end{align}
\end{figure*}
\setcounter{equation}{\value{TempEqCnt}}

\begin{proposition}
\label{prop:zeroDua} If (P2) is feasible, and the selling price $\beta^t$
does not exceed the purchase price $\alpha^t$ for any $t \in \mathcal{T}$,
then there is no duality gap.
\end{proposition}
\begin{IEEEproof}
Due to the strong duality theorem for the optimization problems with
linear constraints~(cf.~\cite[Prop.~5.2.1]{Bertsekas99}), it
suffices to show that the cost function is convex over the entire
space and its optimal value is finite. First, using the same
argument, convexity of $G(\{\tilde{P}_R^t\})$ in $\{\tilde{P}_R^t\}$
is immediate under the transaction price condition. The finiteness
of the optimal value is guaranteed by the fact that the continuous convex
cost~\eqref{P2ObjFunc} is minimized over a nonempty compact set specified by
\eqref{Plimits}--\eqref{Balance}, and~\eqref{VarChange}.
\end{IEEEproof}

The strong duality asserted by Proposition~\ref{prop:zeroDua}
motivates the use of Lagrangian relaxation techniques in order to
solve the scheduling problem. Moreover, problem (P2) is clearly
separable, meaning that its cost and constraints are sums of terms,
with each term dependent on different optimization variables.
%namely $\{P_{G_m}^t\}$, $\{P_{D_n}^t\}$, $\{P_{E_q}^t\}$, $\{B_j^t\}$, $\{P_R^t\}$ and $\{\tilde{P}_R^t\}$.
The features of strong duality and separability
imply that Lagrangian relaxation and dual decomposition are
applicable to yield a decentralized algorithm;
see also related techniques in power systems~\cite{Conejo} and communication networks~\cite{Palomar06,ChiangLCD07}. %Low99
Coordinated by dual variables, the dual approach decomposes the
original problem into several separate subproblems that can be
solved by the LCs in parallel. The development of the distributed
algorithm is undertaken next.

\subsection{Dual Decomposition}

Constraints~\eqref{SpinRes},~\eqref{Balance}, and~\eqref{VarChange}
couple variables across generators, loads, and the RES. Let
$\mathbf{z}$ collect dual variables $\{\mu^t\}$, $\{\lambda^t\}$,
and $\{\nu^t\}$, which denote the corresponding Lagrange multipliers.
%associated with~\eqref{SpinRes},~\eqref{Balance}, and~\eqref{VarChange}, respectively.
Keeping the remaining constraints implicit, the partial Lagrangian is given by
\begin{align}
\hspace{-1mm} &\mathcal{L}(\mathbf{x},\mathbf{z})
=\sum^{T}_{t=1}\Bigg(\sum^{M}_{m=1}C_m^t(P_{G_m}^t)-\sum^{N}_{n=1}U_{D_n}^t(P_{D_n}^t)\nonumber \\
&-\sum^{Q}_{q=1}U_{E_q}^t(P_{E_q}^t)+\sum^{J}_{j=1}H_j^t(B_j^t)\Bigg)+G(\{\tilde{P}_R^t\}) \nonumber \\
&+ \sum^{T}_{t=1}\Bigg\{\mu^t\Bigg(\mathsf{SR}^t-\sum^{M}_{m=1}(P_{G_m}^{\max}-P_{G_m}^{t})\Bigg) \nonumber \\
&-\lambda^t\Bigg(\sum^{M}_{m=1}P_{G_m}^t+P_R^t-\sum^{N}_{n=1}P_{D_n}^t-\sum^{Q}_{q=1}P_{E_q}^t-L^t\Bigg)\nonumber \\
&-\nu^t\Bigg(\tilde{P}_R^t-P_R^t-\sum_{j=1}^JP_{B_j}^t\Bigg)
\Bigg\}.
\end{align}
Then, the dual function can be written as
\begin{align*}
\label{dualFun}
\mathcal{D}(\mathbf{z}) =
&\mathop{\min}_{\mathbf{x}} \mathcal{L}(\mathbf{x},\mathbf{z}) \nonumber \\
&\st~\eqref{Plimits}-\eqref{RampDown},~\eqref{Dlimits}-\eqref{Rlimits}
\end{align*}
and the dual problem is given by
\begin{subequations}
\begin{align}
\mathop{\max} ~~ &\mathcal{D}(\{\mu^t\},\{\lambda^t\},\{\nu^t\}) \\
\st \quad &\mu^t \ge 0, \lambda^t,\nu^t\in \mathbb{R}, ~t \in \cT.
\end{align}
\end{subequations}

The subgradient method will be employed to obtain the optimal multipliers and power schedules.
The iterative process is described next, followed by its distributed implementation.

\subsubsection{Subgradient Iterations}

The subgradient method amounts to running the
recursions~\cite[Sec.~6.3]{Bertsekas09}
\begin{subequations}
\label{LMupd}
\begin{align}
\mu^t(k+1)&=[\mu^t(k)+a g_{\mu^t}(k)]^{+}  \\
\lambda^t(k+1)&=\lambda^t(k)+a g_{\lambda^t}(k)\\
\nu^t(k+1)&=\nu^t(k)+a g_{\nu^t}(k)
\end{align}
\end{subequations}
where $k$ is the iteration index; $a>0$ is a constant stepsize;
while $g_{\mu^t}(k)$,~$g_{\lambda^t}(k)$, and $g_{\nu^t}(k)$ denote
the subgradients of the dual function with respect to
$\mu^t(k)$,~$\lambda^t(k)$, and $\nu^t(k)$, respectively. These
subgradients can be expressed in the following simple forms
\begin{subequations}
\label{Sg}
\begin{align}
g_{\mu^t}(k) &= \mathsf{SR}^t-\sum^{M}_{m=1}(P_{G_m}^{\max}-P_{G_m}^{t}(k))  \\
g_{\lambda^t}(k) &= L^t+\sum^{N}_{n=1}P_{D_n}^t(k)+\sum^{Q}_{q=1}P_{E_q}^t(k)
\nonumber \\ &\hspace{0.4cm}
-\sum^{M}_{m=1}P_{G_m}^t(k) -P_R^t(k)\\
g_{\nu^t}(k) &= P_R^t(k)+\sum_{j=1}^JP_{B_j}^t(k)-\tilde{P}_R^t(k)
\end{align}
\end{subequations}
where $P_{G_m}^{t}(k)$, $P_{D_n}^t(k)$, $P_{E_q}^t(k)$,
$P_{B_j}^t(k)$, $P_R^t(k)$, and $\tilde{P}_R^t(k)$ are given by
\eqref{PGsubp}--\eqref{PRsubp}.
%{\footnotesize
%\begin{align}
%&\{P_{G_m}^{t}(k)\}_{t=1}^T \in \argmin\limits_{\substack{\{P_{G_m}^t\} \\ \st~\eqref{Plimits}-\eqref{RampDown}}}
%\Bigg\{\sum^{T}_{t=1}\bigg(C_m^t(P_{G_m}^t)+\big(\mu^t(k)-\lambda^t(k)\big)P_{G_m}^t\bigg)\Bigg\} \label{PGsubp}\\
%&\{P_{D_n}^t(k)\}_{t=1}^T  \in \argmin\limits_{\substack{\{P_{D_n}^t\} \\ \st~\eqref{Dlimits}}}
%\Bigg\{\sum^{T}_{t=1}\bigg(\lambda^t(k)P_{D_n}^t-U_{D_n}^t(P_{D_n}^t)\bigg)\Bigg\} \label{PDsubp} \\
%&\{P_{E_q}^t(k)\}_{t=1}^T \in \argmin\limits_{\substack{\{P_{E_q}^t\} \\ \st~\eqref{Elimits}-\eqref{Esumlimits}}}
%\Bigg\{\sum^{T}_{t=1}\bigg(\lambda^t(k)P_{E_q}^t-U_{E_q}^t(P_{E_q}^t)\bigg)\Bigg\} \label{PEsubp} \\
%&\{P_{B_j}^t(k)\}_{t=1}^T \in \argmin\limits_{\substack{\{P_{B_j}^t,B_j^t\} \\ \st~\eqref{Blimits}-\eqref{FinalSoC}}}
%\Bigg\{\sum^{T}_{t=1}\bigg(\nu^t(k)P_{B_j}^t+H_j^t(B_j^t)\bigg)\Bigg\} \label{PBsubp} \\
%&\{P_R^t(k),\tilde{P}_R^t(k)\}_{t=1}^T \in \argmin\limits_{\substack{\{P_{R}^t~\st~\eqref{Rlimits}, \\ \tilde{P}_R^t \in \mathbb{R}^{T}\}}}
%\Bigg\{\sum^{T}_{t=1}\bigg(\big(\nu^t(k)-\lambda^t(k)\big)P_R^t\bigg) \nonumber \\
%&\hspace{3cm} +G(\{\tilde{P}_R^t\})-\sum^{T}_{t=1}\nu^t(k)\tilde{P}_R^t\Bigg\}\label{PRsubp}
%\end{align}
%}

Iterations are initialized with arbitrary $\lambda^t(0), \nu^t(0)
\in \mathbb{R}$, and $\mu^t(0)\geq 0$. The iterates are guaranteed
to converge to a neighborhood of the optimal
multipliers~\cite[Sec.~6.3]{Bertsekas09}. The size of the
neighborhood is proportional to the stepsize, and can therefore be
controlled by the stepsize.

When the primal objective is \emph{not} strictly convex, a primal
averaging procedure is necessary to obtain the optimal power
schedules, which are then given by
\setcounter{equation}{20}
\begin{align}
\label{RunAve}
\bar{\mathbf{x}}(k)=\frac{1}{k} \sum_{j=0}^{k-1} \mathbf{x}(j) = \frac{1}{k}\mathbf{x}(k-1)+\frac{k-1}{k}\bar{\mathbf{x}}(k-1).
\end{align}
The running averages can be recursively computed as
in~\eqref{RunAve}, and are also guaranteed to converge to a
neighborhood of the optimal solution~\cite{NedicOzdaglar-PrimSol-JOpt}.
Note that other convergence-guaranteed stepsize rules and primal averaging methods
can also be utilized; see~\cite{Gatsis-TSG} for detailed discussions.

%Note that other convergence-guaranteed stepsize rules (e.g., diminishing stepsize)
%and primal averaging (weights proportional to the stepsize) can also
%be utilized; see~\cite{Gatsis-TSG} for detailed discussions.

\subsubsection{Distributed Implementation}

The form of the subgradient iterations easily lends itself to a
distributed implementation utilizing the control and communication
capabilities of a typical microgrid.

\begin{figure}[ht]
\centering
\includegraphics[scale=0.75]{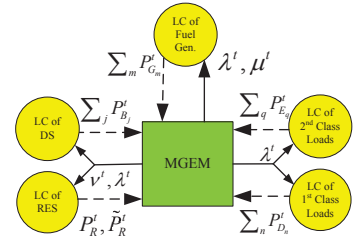}
\caption{Decomposition and message exchange.} \label{fig:decomp}
\vspace*{-0.1cm}
\end{figure}

Specifically, the MGEM maintains and updates the Lagrange
multipliers via~\eqref{LMupd}. The LCs of conventional generation,
dispatchable loads, storage units, and RES solve
subproblems~\eqref{PGsubp}--\eqref{PRsubp}, respectively. These subproblems can be solved if
the MGEM sends the current multiplier iterates $\mu^t(k)$,
$\lambda^t(k)$, and $\nu^t(k)$  to the LCs. The LCs send back to the
MGEM the quantities $\sum^{M}_{m=1}P_{G_m}^{t}(k)$,
$\sum^{N}_{n=1}P_{D_n}^t(k)$, $\sum^{Q}_{q=1}P_{E_q}^t(k)$,
$\sum_{j=1}^JP_{B_j}^t(k)$, $P_R^t(k)$, and $\tilde{P}_R^t(k)$ which
are in turn used to form the subgradients according to~\eqref{Sg}.
The distributed algorithm using dual decomposition is tabulated as
Algorithm~\ref{algo:Distri}, and the interactive process of message
passing is illustrated in Fig.~\ref{fig:decomp}.

%%%%%%%%%%%%%%%%%%%%%%%%%%%%%%%%%%%%%%%%%%%%%%%%%%%%
\begin{algorithm}[t]
\caption{Distributed Energy Management}
\label{algo:Distri}
\begin{algorithmic}[1]
\State Initialize Lagrange multipliers $\lambda^t = \mu^t = \nu^t = 0$
\Repeat \quad ($k = 0,1,2,\ldots$)
\For{$t=1,2,\dots,T$}
    \State Broadcast $\lambda^t(k)$, $\mu^t(k)$, and $\nu^t(k)$ to LCs of convectional generators, controllable loads, storage units, and RES facilities
    \State Update power scheduling $P_{G_m}^t(k)$, $P_{D_n}^t(k)$, $P_{E_q}^t(k)$, $P_{B_j}^t(k)$, $P_{R}^t(k)$, and $\tilde{P}_R^t(k)$ by solving~\eqref{PGsubp}--\eqref{PRsubp}
    \State Update $\lambda^t(k)$,~$\mu^t(k)$, and $\nu^t(k)$ via~\eqref{LMupd}
\EndFor
\State Running averages of primal variables via~\eqref{RunAve}
\Until Convergence
\end{algorithmic}
\end{algorithm}
%%%%%%%%%%%%%%%%%%%%%%%%%%%%%%%%%%%%%%%%%%%%%%%%%%%%

%%%%%%%%%%%%%%%%%%%%%%%%%%%%%%%%%%%%%%%%%%%%%%%%%%%%
\begin{algorithm}[ht]
\caption{Enumerate all the vertices of a polytope $\cA$}
\label{algo:ListV}
\begin{algorithmic}[1]
\State Initialize vertex set $\cV = \emptyset$ \State Generate set
$\tilde{\cA}:=\{\tilde{\mathbf{a}} \in \mathbb{R}^n|\tilde{a}_i
=\underline{a}_i~\textrm{or}~\overline{a}_i,~i=1,\ldots,n\}$; check
the feasibility of all the points in set $\tilde{\cA}$, i.e., if
$a^{\min}\le \mathbf{1}^{\prime}\tilde{\mathbf{a}} \le a^{\max}\}$,
then $\cV = \cV \cup \{\tilde{\mathbf{a}}\}$ \State Generate set
$\hat{\cA}:=\{\hat{\mathbf{a}} \in \mathbb{R}^n| \hat{a}_i =
a^{\min}-\sum_{j \neq i}\hat{a}_j~\textrm{or}~a^{\max}-\sum_{j \neq
i}\hat{a}_j,~\hat{a}_j
=\underline{a}_j~\textrm{or}~\overline{a}_j,~i, j=1,\ldots,n, j \neq
i\}$; check the feasibility of all the points in set $\hat{\cA}$,
i.e., if $\underline{\mathbf{a}}\preceq \hat{\mathbf{a}} \preceq
\overline{\mathbf{a}}$, then $\cV = \cV \cup \{\hat{\mathbf{a}}\}$
\end{algorithmic}
\end{algorithm}
%%%%%%%%%%%%%%%%%%%%%%%%%%%%%%%%%%%%%%%%%%%%%%%%%%%%

%%%%%%%%%%%%%%%%%%%%%%%%%%%%%%%%%%%%%%%%%%%%%%%%%%%%
\begin{algorithm}[ht]
\caption{Enumerate all the vertices of a polytope $\cB$}
\label{algo:CombV}
\begin{algorithmic}[1]
\For{$i=1,2,\ldots,S$}
  \State Obtain vertex set $\cV_s$ by applying Algorithm~\ref{algo:ListV} to $\cB_s$
\EndFor
\State Generate vertices $\mathbf{b}^{\tV}$ for $\cB$ by concatenating all the individual vertices $\mathbf{b}_s$ as
$\mathbf{b}^{\tV} =[(\mathbf{b}_1^{\tV})^{\prime},\ldots,(\mathbf{b}_S^{\tV})^{\prime}]^{\prime},~\mathbf{b}_s \in \cV_s$
\end{algorithmic}
\end{algorithm}
%%%%%%%%%%%%%%%%%%%%%%%%%%%%%%%%%%%%%%%%%%%%%%%%%%%%

\subsection{Solving the LC Subproblems}

This subsection shows how to solve each
subproblem~\eqref{PGsubp}--\eqref{PRsubp}. Specifically,
$C_m^t(\cdot)$, $-U_{D_n}^t(\cdot)$, $-U_{E_q}^t(\cdot)$, and
$H_j^t(\cdot)$ are chosen either convex piece-wise linear or smooth
convex quadratic. Correspondingly, the first four subproblems~\eqref{PGsubp}--\eqref{PBsubp} are
essentially linear programs (LPs) or quadratic programs (QPs), which can be solved efficiently.
Therefore, the main focus is on solving~\eqref{PRsubp}.

The optimal solution of $P_R^t(k)$ in~\eqref{PRsubp} is easy to obtain as
\begin{align}
\label{PR(k)}
P_R^t(k)  =
\left\{\begin{array}{cc}
P_{R}^{\min},  &\mbox{if}~\nu^t(k) \ge \lambda^t(k) \\
P_{R}^{\max},  &\mbox{if}~\nu^t(k) < \lambda^t(k). \\
\end{array}\right.
\end{align}
However, due to the absolute value operator and the maximization over $\mathbf{w}$
in the definition of $G(\{\tilde{P}_R^t\})$, subproblem~\eqref{PRsubp} is a convex
nondifferentiable problem in $\{\tilde{P}_R^t\}$, which can be
challenging to solve. As a state-of-the-art technique for convex
nondifferentiable optimization
problems~\cite[Ch.~6]{Bertsekas09}, the bundle
method is employed to obtain $\{\tilde{P}_R^t(k)\}$.

Upon defining
\begin{align}
\tilde{G}(\{\tilde{P}_R^t\}) := G(\{\tilde{P}_R^t\})-\sum^{T}_{t=1}\nu^t(k)\tilde{P}_R^t
\end{align}
the subgradient of $\tilde{G}(\{\tilde{P}_R^t\})$ with respect to $\tilde{P}_R^t$ needed for the bundle method
can be obtained by the generalization of Danskin's Theorem~\cite[Sec.~6.3]{Bertsekas09} as
\begin{align}
\label{subgrad}
\partial\tilde{G}(\{\tilde{P}_R^t\})  =
\left\{\begin{array}{cc}
\alpha^t-\nu^t(k),  &\mbox{if}~\tilde{P}_R^t \ge \sum\limits_{i=1}^{I}(W_i^t)^{*} \\
\beta^t-\nu^t(k),  &\mbox{if}~\tilde{P}_R^t <\sum\limits_{i=1}^{I}(W_i^t)^{*}
\end{array}\right.
\end{align}
where for given $\{\tilde{P}_R^t\}$ it holds that
{\small
\begin{align}
\label{prob:Wopt}
\mathbf{w}^{*} \in
\argmax\limits_{\mathbf{w} \in \cW}
\left\{\sum_{t=1}^T\left(\delta^t\left|\tilde{P}_R^t-\sum_{i=1}^{I}W_i^t\right|
+\gamma^t\left(\tilde{P}_R^t-\sum_{i=1}^{I}W_i^t\right)\right)\right\}.
\end{align}}
With $\mathbf{p}:= [\tilde{P}_R^1,\ldots,\tilde{P}_R^T]$, the bundle
method generates
%$\mathbf{p}^{\min}:= P_R^{\min}\cdot\mathbf{1}$, $\mathbf{p}^{\max}:= P_R^{\max}\cdot\mathbf{1}$,where $\mathbf{1}$ is the all-ones vector.
a sequence $\{\mathbf{p}_{\ell}\}$ with guaranteed convergence to
the optimal $\{\tilde{P}_R^t(k)\}$; see e.g.,~\cite{Kiwiel00},~\cite[Ch.~6]{Bertsekas09}.
The iterate $\mathbf{p}_{\ell+1}$ is obtained by minimizing a polyhedral approximation
%(piece-wise linear approximation)
of $\tilde{G}(\mathbf{p})$ with a quadratic proximal regularization
as follows
\begin{align}
\label{prob:bundle}
\mathbf{p}_{\ell+1} :=\argmin\limits_{\mathbf{p} \in \mathbb{R}^T} \left\{\hat{G}_{\ell}(\mathbf{p})+\frac{\rho_{\ell}}{2}\|\mathbf{p}-\mathbf{y}_{\ell}\|^2\right\}
\end{align}
where
$\hat{G}_{\ell}(\mathbf{p}):=\max\{\tilde{G}(\mathbf{p}_{0})+\mathbf{g}_{0}^{\prime}(\mathbf{p}-\mathbf{p}_{0}),
\ldots,\tilde{G}(\mathbf{p}_{\ell})+\mathbf{g}_{\ell}^{\prime}(\mathbf{p}-\mathbf{p}_{\ell})\}$;
$\mathbf{g}_{\ell}$ is the subgradient of $\tilde{G}(\mathbf{p})$
evaluated at the point $\mathbf{p}=\mathbf{p}_{\ell}$, which is
calculated according to~\eqref{subgrad}; proximity weight
$\rho_{\ell}$ is to control stability of the iterates; and the
proximal center $\mathbf{y}_{\ell}$ is updated according to a query
for descent
\begin{align}
\mathbf{y}_{\ell+1}  =
\left\{\begin{array}{cc}
\mathbf{p}_{\ell+1},  &\mbox{if}~\tilde{G}(\mathbf{y}_{\ell})-\tilde{G}(\mathbf{p}_{\ell+1}) \ge \theta\eta_{\ell} \\
\mathbf{y}_{\ell},  &\mbox{otherwise}
\end{array}\right.
\end{align}
where $\eta_{\ell} = \tilde{G}(\mathbf{y}_{\ell})-\left(\hat{G}_{\ell}(\mathbf{p}_{\ell+1})+
\frac{\rho_{\ell}}{2}\|\mathbf{p}_{\ell+1}-\mathbf{y}_{\ell}\|^2\right)$, $\theta \in (0,1)$.

It is worth mentioning that~\eqref{prob:bundle} is essentially a QP over
a simplex in the dual space, which is efficiently solvable by practical optimization algorithms.
The corresponding transformation is shown in Appendix~\ref{Appe:Bundle}
for the interested readers.

Algorithms for solving~\eqref{prob:Wopt} depend on the form of the
uncertainty set $\cW$, and are elaborated next.

\subsection{Vertex Enumerating Algorithms}
In order to obtain $\mathbf{w}^{*}$, the convex nondifferentiable
function in~\eqref{prob:Wopt} should be maximized over $\cW$. This
is generally an NP-hard convex maximization problem.
%meaning the global optimum $\{W_{*}^t\}$ can not be obtained in polynomial time.
However, for the specific problem here, the special structure of the problem can be utilized
to obtain a computationally efficient approach.

Specifically, the global solution is attained at the extreme points of the
polytope~\cite[Sec.~2.4]{Bertsekas09}. Therefore, the objective in~\eqref{prob:Wopt} can be
evaluated at all vertices of $\cW$ to obtain the global solution.
Since there are only finitely many vertices,~\eqref{prob:Wopt} can
be solved in a \emph{finite} number of steps.

For the polytopes $\cW$ with special structure
[cf.~\eqref{perW},~\eqref{jointW}], characterizations of vertices
are established in Propositions~\ref{prop:Vchara}
and~\ref{prop:Vcomb}. Capitalizing on these propositions, vertex
enumerating procedures are designed consequently, and are tabulated
as Algorithms~\ref{algo:ListV} and~\ref{algo:CombV}.
\begin{proposition}
\label{prop:Vchara} For a polytope $\cA := \{\mathbf{a} \in
\mathbb{R}^n|\underline{\mathbf{a}}\preceq \mathbf{a} \preceq
\overline{\mathbf{a}}, a^{\min}\le  \mathbf{1}^{\prime}\mathbf{a}
\le a^{\max}\}$, $\mathbf{a}^{\tV} \in \cA$ is a vertex (extreme
point) of $\cA$ if and only if it has one of the following forms: i)
${a}^{\tV}_i =\underline{a}_i~\textrm{or}~\overline{a}_i$ for
$i=1,\ldots,n$; or ii) ${a}^{\tV}_i = a^{\min}-\sum_{j \neq
i}{a}^{\tV}_j~\textrm{or}~a^{\max}-\sum_{j \neq i}{a}^{\tV}_j,
~{a}^{\tV}_j =\underline{a}_j~\textrm{or}~\overline{a}_j$, for $i,
j=1,\ldots,n, j \neq i$.
\end{proposition}
\begin{IEEEproof}
See Appendix~\ref{Appe:Vchara}.
\end{IEEEproof}
Essentially, Proposition~\ref{prop:Vchara} verifies the geometric
characterization of vertices. Since $\cW$ is the part of a
hyperrectangle (orthotope) between two parallel hyperplanes, its
vertices can only either be the hyperrectangle's vertices which are
not cut away, or, the vertices of the intersections of the
hyperrectangle and the hyperplanes, which must appear in some edges
of the hyperrectangle.

\begin{table}[t]
\centering
\caption{Generating capacities, ramping limits, and cost coefficients.
The units of $a_m$ and $b_m$ are \$/(kWh)$^{2}$ and \$/kWh, respectively.}\label{tab:generator}
    \begin{tabular}{  c || c | c | c | c | c }
    \hline
Unit &$P_{G_m}^{\min}$ &$P_{G_m}^{\max}$  &$R_{m,\text{up(down)}}$  &$a_m$   &$b_m$       \\ \hline
1    & 10                     & 50             & 30             &0.006                  & 0.5                      \\
2    & 8                      & 45             & 25             &0.003                  & 0.25                    \\
3    & 15                     & 70             & 40              &0.004                  & 0.3                         \\
    \hline
    \end{tabular}
\end{table}

\begin{table}[t]
\centering
\caption{Class-1 dispatchable loads parameters.
The units of $c_n$ and $d_n$ are \$/(kWh)$^{2}$ and \$/kWh, respectively.}\label{tab:loadK1}
    \begin{tabular}{ c || c | c | c | c | c | c }
    \hline
                       &Load 1   &Load 2    &Load 3    &Load 4       &Load 5    &Load 6   \\  \hline \hline
    $P_{D_n}^{\min}$                         & 0.5         & 4        & 2        & 5.5       & 1        & 7     \\
    $P_{D_n}^{\max}$                        & 10         & 16       & 15        & 20       & 27        & 32   \\ \hline
    $c_n$                       & -0.002      & -0.0017    & -0.003     & -0.0024     & -0.0015      & -0.0037     \\
    $d_n$                       & 0.2         & 0.17     &0.3        &0.24       &0.15        &0.37   \\
\hline
    \end{tabular}
\end{table}

\begin{table}[t]
\centering
\caption{Class-2 dispatchable loads parameters}\label{tab:loadK2}
    \begin{tabular}{ c || c | c | c | c }
    \hline
                          &Load 1        &Load 2        &Load 3       &Load 4    \\  \hline \hline
    $P_{E_q}^{\min}$                &0           &0               &0           &0          \\
    $P_{E_q}^{\max}$                &1.2        &1.55          &1.3         &1.7        \\
    $E_q^{\max}$                    &5          &5.5             &4           &8        \\
    $S_q$                           &6PM        &7PM            &6PM          &6PM          \\
    $T_q$                           &12AM       &11PM           &12AM         &12AM        \\
\hline
    \end{tabular}
\end{table}

\begin{table}[t]
\renewcommand{\arraystretch}{1.2}
\centering
\caption{Limits of forecasted wind power}\label{tab:wind}
    \begin{tabular}{  c || c | c | c | c | c | c | c | c }
    \hline
    \text{Slot}             & 1    &  2   &  3    &   4  &   5  &  6    &   7   & 8 \\ \hline \hline
    $\underline{W}_1^t$     &2.47&2.27&2.18&1.97&2.28&2.66&3.1&3.38 \\
    $\overline{W}_1^t$      &24.7&22.7&21.8&19.7&22.8&26.6&31&33.8 \\
    $\underline{W}_2^t$     &2.57&1.88&2.16&1.56&1.95&3.07&3.44&3.11  \\
    $\overline{W}_2^t$      &25.7&18.8&21.6&15.6&19.5&30.7&34.4&31.1 \\
    \hline
    \end{tabular}
\end{table}

\begin{table}[t]
\centering
\caption{Fixed loads demand and transaction prices.
The units of $\alpha^t$ and $\beta^t$ are \textcent/kWh.}\label{tab:price}
   %\begin{tabular}{  c || c | c | c | c | c | c | c | c }
    \begin{tabular}{p{11mm}||p{5mm}|p{5mm}|p{5mm}|p{5mm}|p{7mm}|p{5mm}|p{5mm}|p{5mm}}
    \hline
    \text{Slot}             & 1    &  2   &  3    &   4  &   5  &  6    &   7   & 8 \\ \hline \hline
    $L^t$                   &57.8&58.4&64&65.1&61.5&58.8&55.5&51 \\ \hline
    (Case A)              &     &     &      &     &     &      &      & \\
    $\alpha^t$              &2.01&2.2&3.62&6.6&5.83&3.99&2.53&2.34\\
    $\beta^t$    &1.81&1.98&3.26&5.94&5.25&3.59&2.28&2.11 \\ \hline
    (Case B)             &     &     &      &     &     &      &      & \\
    $\alpha^t$             &40.2&44&72.4&132&116.6&79.8&50.6&46.8\\
    $\beta^t$     &36.18&39.6&65.16&118.8&104.94&71.82&45.54&42.12 \\
    \hline
    \end{tabular}
\end{table}

Next, the vertex characterization of a polytope in a Cartesian product
formed by many lower-dimensional polytopes like $\cA$ is established, which
is needed for the uncertainty set~\eqref{perW}.

\begin{proposition}
\label{prop:Vcomb} Assume $\mathbf{b}\in \mathbb{R}^n$ is divided
into $S$ consecutive and non-overlapping blocks as $\mathbf{b} =
[\mathbf{b}_1^{\prime},\ldots,\mathbf{b}_S^{\prime}]^{\prime}$,
where $\mathbf{b}_s \in \mathbb{R}^{n_s}$ and $\sum_{s=1}^{S}n_s=n$.
Consider a polytope $\cB := \{\mathbf{b}\in
\mathbb{R}^n|\underline{\mathbf{b}}\preceq \mathbf{b}\preceq
\overline{\mathbf{b}}, b^{\min}_s \le
\mathbf{1}_{n_s}^{\prime}\mathbf{b}_s \le b^{\max}_s, s =
1,\ldots,S\}$. Then $\mathbf{b}^{\tV}
=[(\mathbf{b}_1^{\tV})^{\prime},\ldots,(\mathbf{b}_S^{\tV})^{\prime}]^{\prime}$
is a vertex of $\cB$ if and only if for $s = 1,\ldots,S$,
$\mathbf{b}_s^{\tV}$ is the vertex of a lower-dimensional polytope
$\cB_s := \{\mathbf{b}_s\in
\mathbb{R}^{n_s}|\underline{\mathbf{b}}_s\preceq \mathbf{b}_s\preceq
\overline{\mathbf{b}}_s, b^{\min}_s \le
\mathbf{1}_{n_s}^{\prime}\mathbf{b}_s \le b^{\max}_s\}$.
\end{proposition}
\begin{IEEEproof}
See Appendix~\ref{Appe:Vcomb}.
\end{IEEEproof}

Algorithms~\ref{algo:ListV} and~\ref{algo:CombV} can be used to
to generate the vertices of uncertainty sets~\eqref{perW} and~\eqref{jointW} as described next.
\begin{itemize}
\item[i)] For uncertainty set~\eqref{perW}, first use Algorithm~\ref{algo:ListV} to obtain the vertices corresponding to each sub-horizon $\cT_{i,s}$
for all the RES facilities. Then, concatenate the obtained vertices to get the ones for each RES facility by Step~4 in Algorithm~\ref{algo:CombV}.
Finally, run this step again to form the vertices of~\eqref{perW} by concatenating the vertices of each $\cW_i$.
\item[ii)] For uncertainty sets~\eqref{jointW}, use Algorithm~\ref{algo:ListV} to obtain the vertices for each sub-horizon $\cT_{s}$.
Note that concatenating step in Algorithm~\ref{algo:CombV} is not needed in this case because
problem~\eqref{prob:Wopt} is decomposable across sub-horizons $\cT_{s}$, $s=1,\ldots,S$, and can be independently solved accordingly.
\end{itemize}

After the detailed description of vertex enumerating procedures for RES uncertainty sets, a discussion
on the complexity of solving~\eqref{prob:Wopt} follows.
\begin{remark}\textit{(Complexity of solving~\eqref{prob:Wopt})}.
Vertex enumeration incurs exponential complexity because the number
of vertices can increase exponentially with the number of variables
and constraints~~\cite[Ch.~2]{BertsimasT97}. However, if the
cardinality of each sub-horizon $\cT_s$ is not very large (e.g.,
when $24$ hours are partitioned into $4$ sub-horizons each
comprising $6$ time slots), then the complexity is affordable. Most
importantly, the vertices of $\cW$ need only be listed once, before
optimization.
\end{remark}

\begin{figure}[t]
\centering
\includegraphics[scale=0.45]{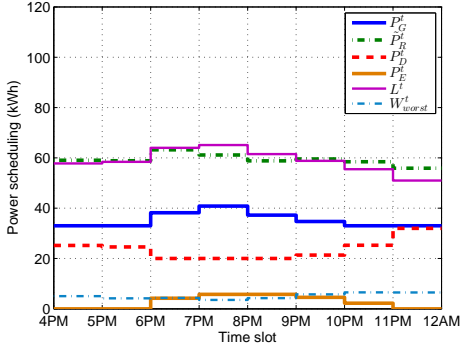}
\caption{Optimal power schedules: Case A.}
\label{fig:powerA}
\end{figure}

\begin{figure}[t]
\centering
\includegraphics[scale=0.45]{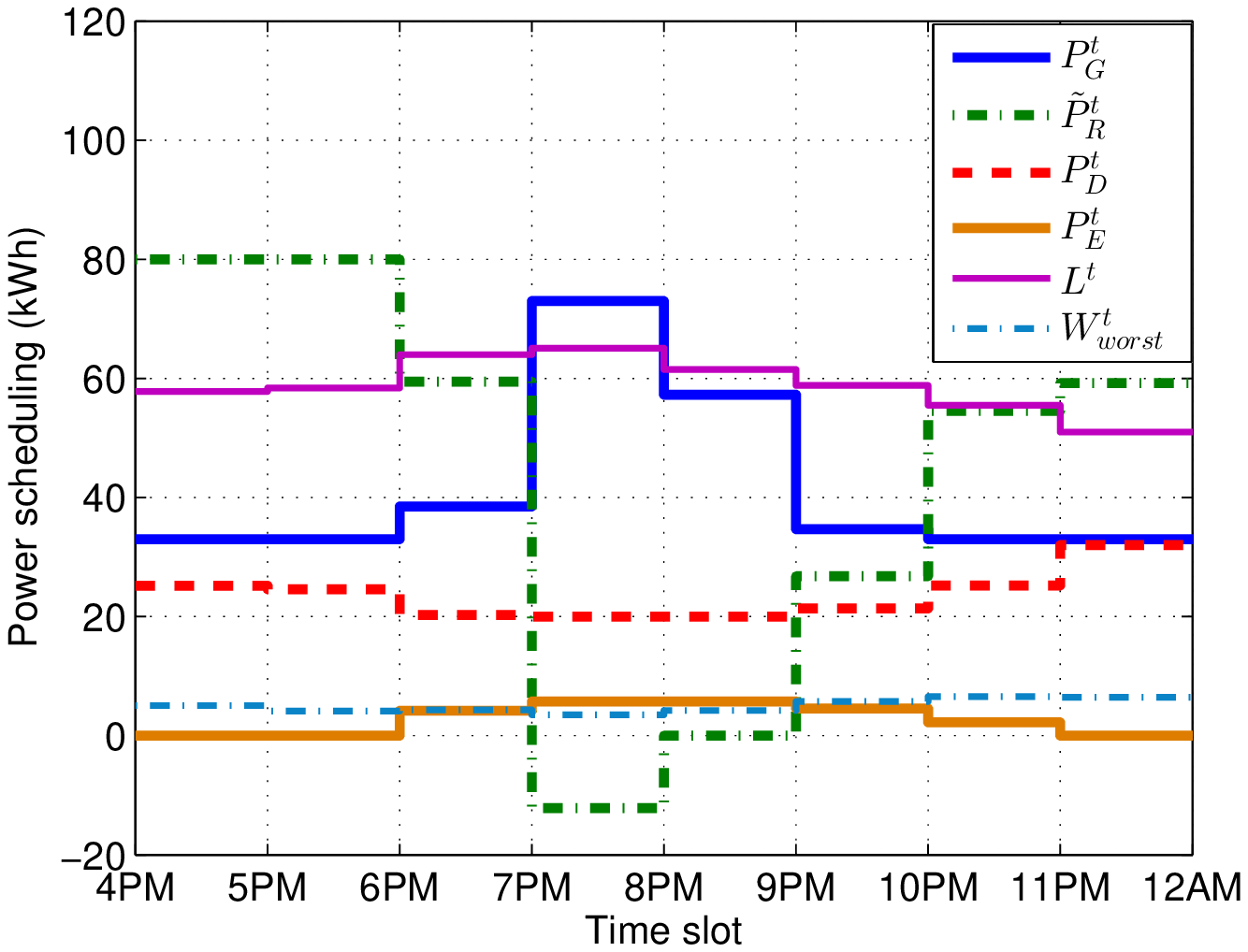}
\caption{Optimal power schedules: Case B.}
\label{fig:powerB}
\end{figure}

\begin{figure}[th]
\centering
\includegraphics[scale=0.45]{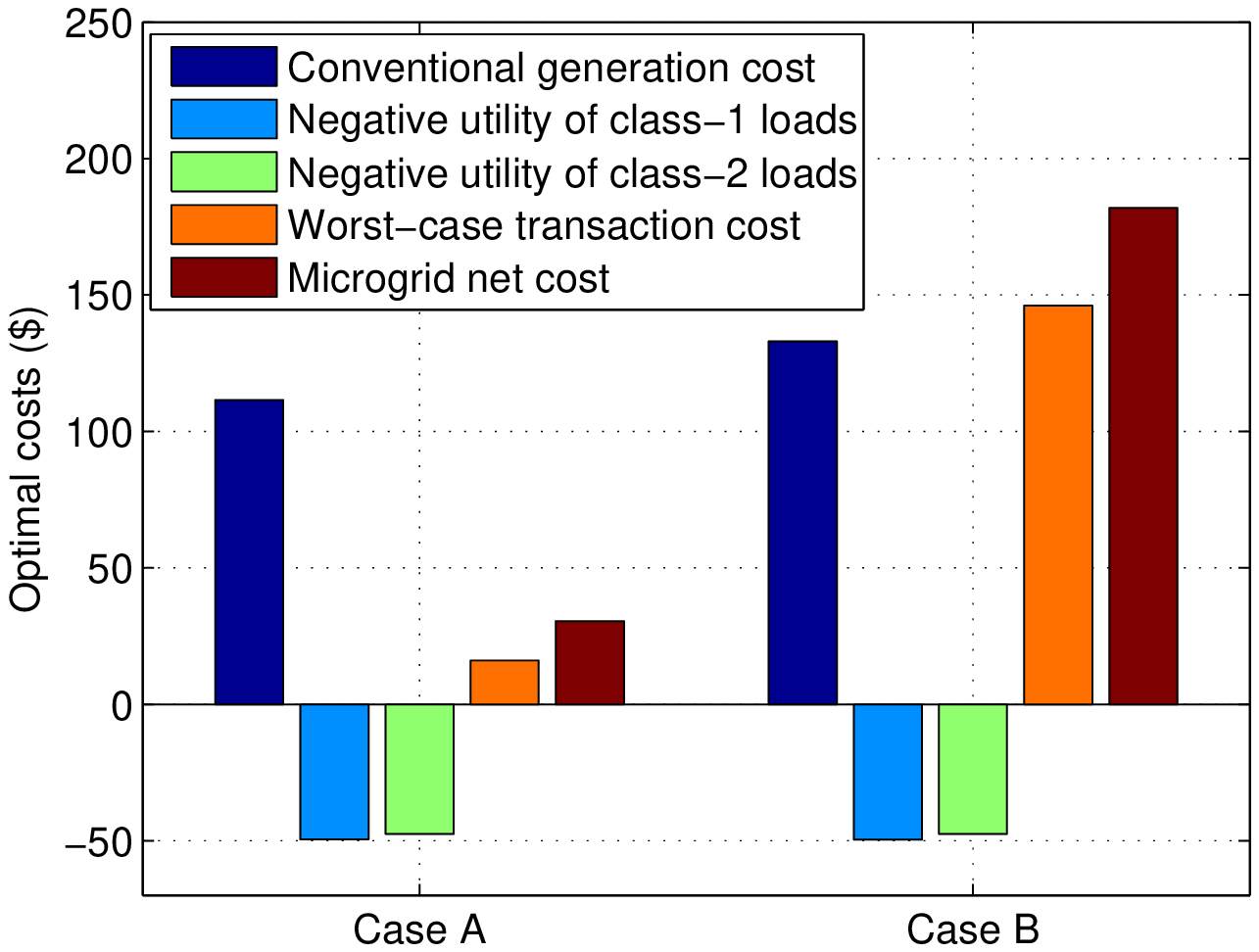}
\caption{Optimal costs: Case A and B.}
\label{fig:costAB}
\end{figure}

\begin{figure}[th]
\centering
\includegraphics[scale=0.45]{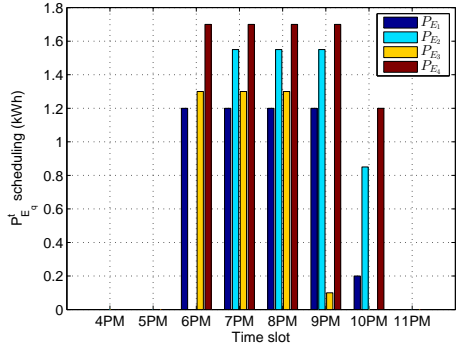}
\caption{Optimal power schedule for $P_{E_q}^t$: Case A.}
\label{fig:PEqA}
\end{figure}

\begin{figure}[th]
\centering
\includegraphics[scale=0.45]{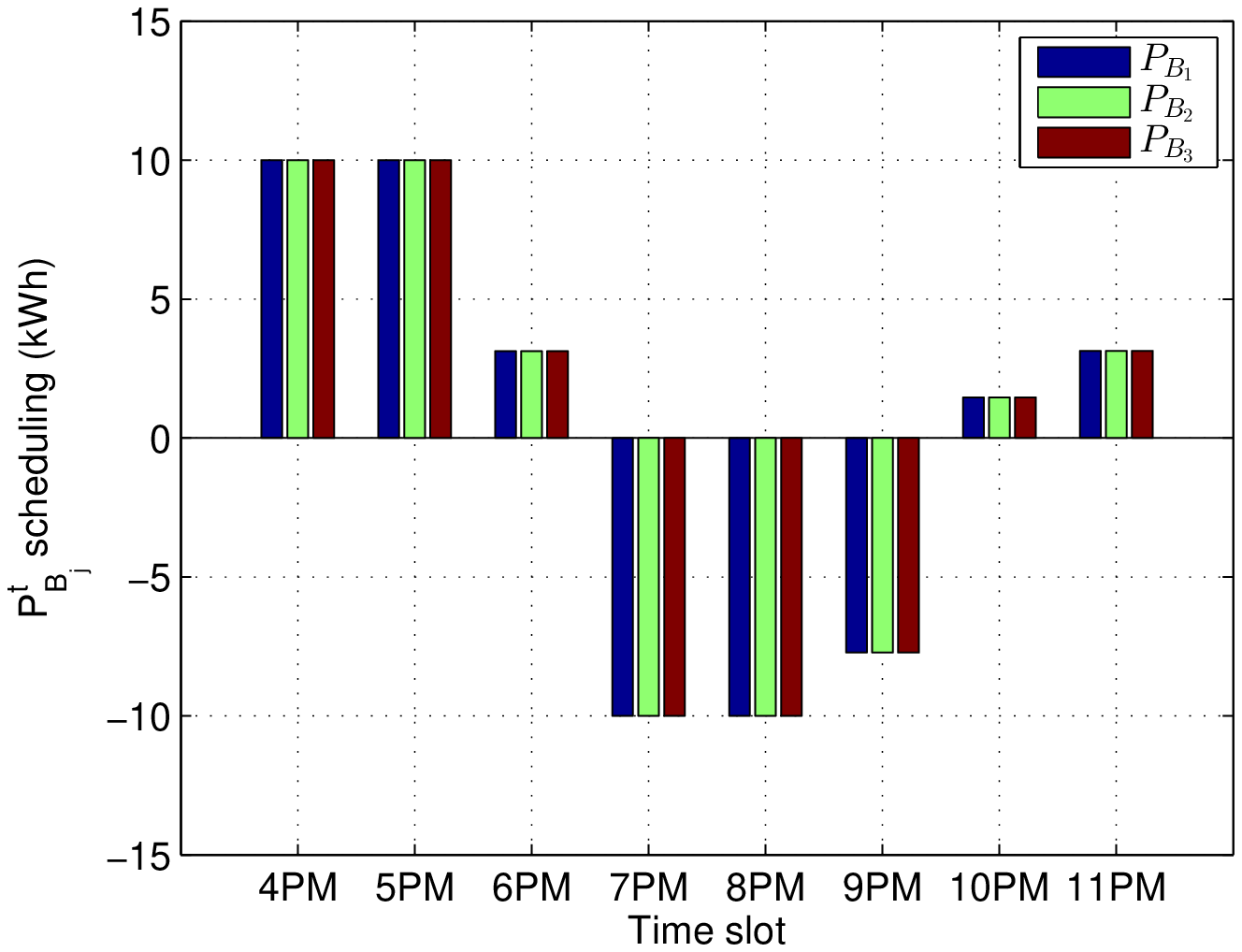}
\caption{Optimal power schedule for $P_{B_j}^t$: Case B.}
\label{fig:PBjB}
\end{figure}

\begin{figure}[th]
\centering
\includegraphics[scale=0.45]{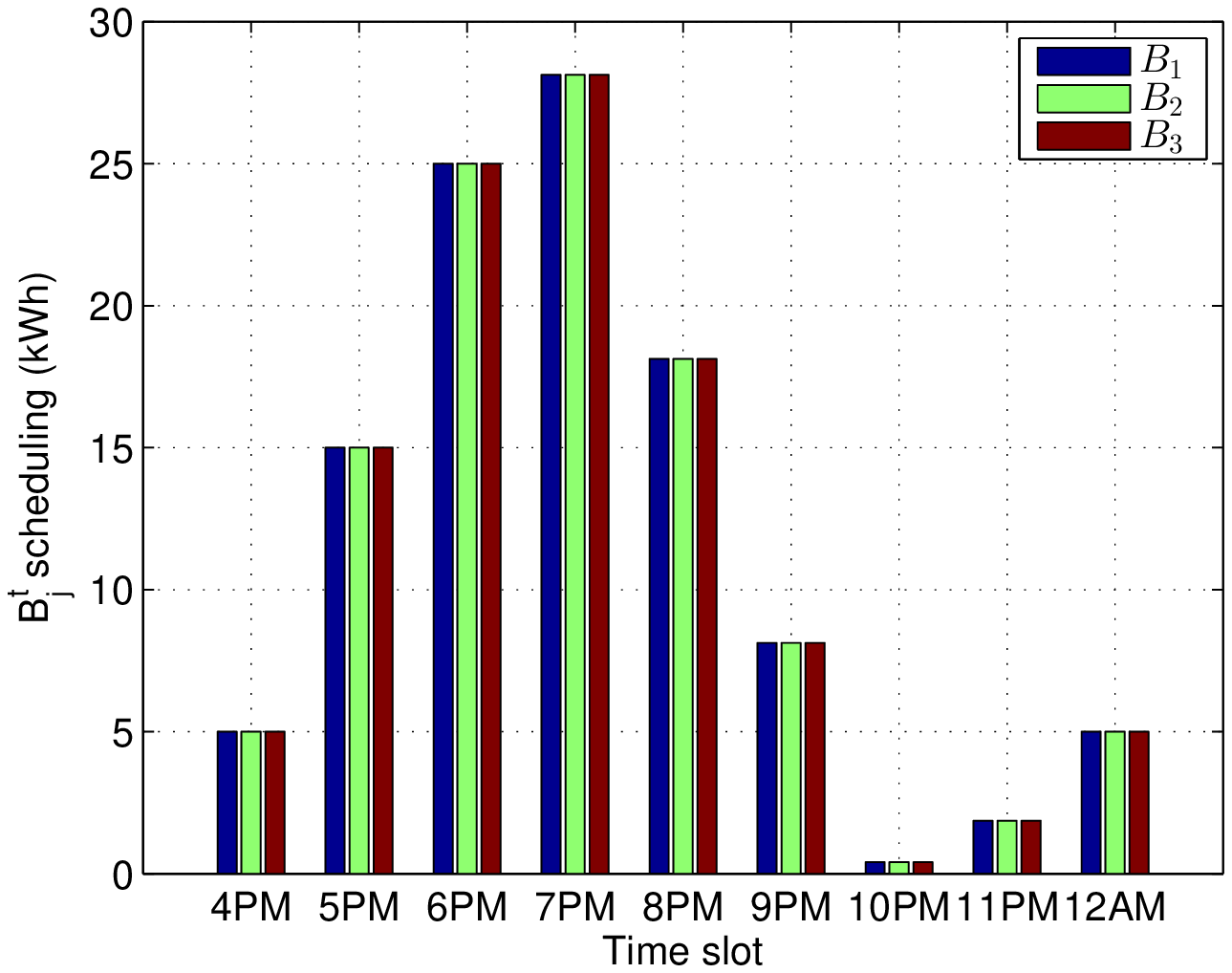}
\caption{Optimal power schedule for $B_j^t$: Case B.}
\label{fig:BjB}
\end{figure}

\begin{figure}[th]
\centering
\includegraphics[scale=0.45]{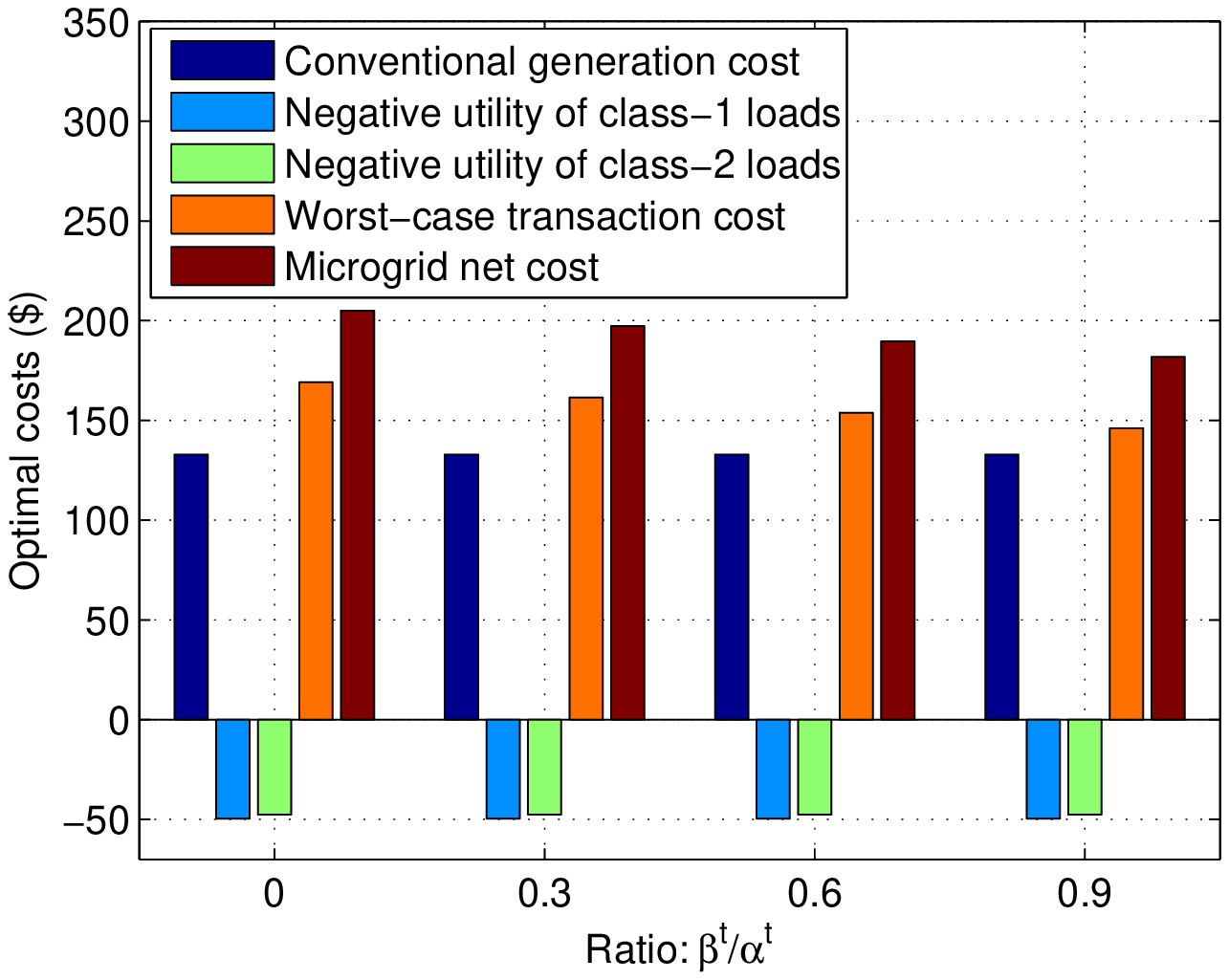}
\caption{Optimal costs: Case B.}
\label{fig:ratiocostB}
\end{figure}

% ----------------------------
\section{Numerical Tests}
\label{sec:Numericalresults}
% ----------------------------

%First, convergence of the Lagrange multiplier~$\{\lambda^t\}$
%corresponding to the balance equation~\eqref{Balance} is confirmed
%for Case A by Fig.~\ref{fig:lambdaA}. It can be seen that
%$\lambda^t$ converges for all $t \in \cT$ within a couple of hundred
%iterations, which verifies the validity of the dual decomposition
%approach using the subgradient projection method. With the running
%averages, convergence of the other dual and primal variables was
%also verified (for Case B too), but is omitted for brevity.
%
%\begin{figure}[ht]
%\centering
%\includegraphics[scale=0.5]{lambdaA}
%\caption{Convergence of $\{\lambda^t\}$.}
%\label{fig:lambdaA}
%\end{figure}

In this section, numerical results are presented to verify the
performance of the robust and distributed energy scheduler.
The Matlab-based modeling package \texttt{CVX}~\cite{cvx} along with the solver \texttt{MOSEK}~\cite{mosek} %\texttt{SeDuMi}~\cite{Sturm99}
are used to specify and solve the proposed robust energy management problem.
The considered microgrid consists of $M=3$ conventional generators,
$N=6$ class-1 dispatchable loads, $Q=4$ class-2
dispatchable loads, $J=3$ storage units, and $I=2$ renewable
energy facilities (wind farms). The time horizon spans $T=8$ hours,
corresponding to the interval $4$PM--$12$AM. The generation costs
$C_m(P_{G_m})=a_mP_{G_m}^2+b_mP_{G_m}$ and the utilities of class-1
elastic loads $U_n(P_{D_n})=c_nP_{D_n}^2+d_nP_{D_n}$ are set
to be quadratic and time-invariant. Generator parameters are given
in Table~\ref{tab:generator}, while $\mathsf{SR}^t =10$kWh. The
relevant parameters of two classes of dispatchable loads are
listed in Tables~\ref{tab:loadK1} and~\ref{tab:loadK2} (see
also~\cite{Gatsis-TSG}). The
utility of class-2 loads is $U_{E_q}^t(P_{E_q}^t) :=
\pi_q^tP_{E_q}^t$ with weights $\pi_q^t=4,3.5,\ldots,1,0.5$ for
$t=4\text{PM},\ldots,11\text{PM}$ and $q\in \cQ$.

Three batteries have capacity $B_j^{\max}=30$kWh (similar to\cite{Guan10}).
The remaining parameters are $P_{B_j}^{\min}=-10$kWh, $P_{B_j}^{\max}=10$kWh, $B_j^0=B_j^{\min}=5$kWh,
and $\eta_j=0.95$, for all $j\in \cJ$.
%The penalty weights for $t=4\text{pm},\ldots,11\text{am}$ are set as
%$\varpi_1^t=0.05,0.1,\ldots,0.35,0.4$ and $\varpi_j^t=0.1$, $j=2,3$.
The battery costs $H_j^t(B_j^t)$ are set to zero.
The joint uncertainty model with $S=1$ is considered for $\cW$ [cf.~\eqref{jointW}],
where $W_1^{\min}=40$kWh, and $W_1^{\max}=360$kWh.
In order to obtain $\underline{W}_i^t$ and $\overline{W}_i^t$ listed in Table~\ref{tab:wind},
MISO day-ahead wind forecast data~\cite{WindData} are rescaled to the order of $1$~kWh to $40$~kWh,
which is a typical wind power generation for a microgrid~\cite{Wu11}.

Similarly, the fixed load $L^t$ in Table~\ref{tab:price} is a rescaled
version of the cleared load provided by MISO's daily
report\cite{MISODaily12}. For the transaction prices, two different
cases are studied as given in Table~\ref{tab:price}, where
$\{\alpha^t\}$ in Case A are real-time prices of the Minnesota hub
in MISO's daily report. To evaluate the effect of high transaction
prices, $\{\alpha^t\}$ in Case B is set as $20$ times of that in Case A.
For both cases, $\beta^t=0.9\alpha^t$, which satisfies the convexity
condition for (P1) given in Proposition~\ref{prop:convex}.

The optimal microgrid power schedules of two cases are shown in
Figs.~\ref{fig:powerA} and~\ref{fig:powerB}. The stairstep curves
include $P_G^t :=\sum_mP_{G_m}^t$, $P_D^t :=\sum_nP_{D_n}^t$, and
$P_E^t :=\sum_qP_{E_q}^t$ denoting the total conventional power
generation, and total elastic demand for classes 1 and
2, respectively, which are the optimal solutions of (P2).
Quantity $W_{\textrm{worst}}^t$ denotes the total worst-case wind energy
at slot $t$, which is the optimal solution of~\eqref{prob:Wopt} with
optimal $\tilde{P}_R^t$.

A common observation from Figs.~\ref{fig:powerA}
and~\ref{fig:powerB} is that the total conventional power generation
$P_G^t$ varies with the same trend across $t$ as the fixed load
demand $L^t$, while the class-1 elastic load exhibits the opposite
trend. Because the conventional generation and the power drawn from
the main grid are limited, the optimal scheduling by solving (P2)
dispatches less power for $P_D^t$ when $L^t$ is large (from $6$PM to
$10$PM), and vice versa. This behavior indeed reflects the load
shifting ability of the proposed design for the microgrid energy
management.

Furthermore, by comparing two cases in Figs.~\ref{fig:powerA}
and~\ref{fig:powerB}, it is interesting to illustrate the effect of
the transaction prices. Remember that the difference between
$\tilde{P}_R^t$ and $W_{\textrm{worst}}^t$ is the shortage power
needed to purchase (if positive) or the surplus power to be sold (if
negative), Figs.~\ref{fig:powerA} shows that the microgrid always
purchases energy from the main grid because $\tilde{P}_R^t$ is more
than $W_{\textrm{worst}}^t$. This is because for Case A, the
purchase price $\alpha^t$ is much lower than the marginal cost of
the conventional generation (cf.~Tables~\ref{tab:generator}
and~\ref{tab:price}). The economic scheduling decision is thus to
reduce conventional generation while purchasing more power to keep
the supply-demand balance. For Case B, since $\alpha^t$ is much
higher than that in Case A, less power should be purchased which is
reflected in the relatively small gap between $\tilde{P}_R^t$ and
$W_{\textrm{worst}}^t$ across time slots. It can also be seen that
$\tilde{P}_R^t$ is smaller than $W_{\textrm{worst}}^t$ from $7$PM to
$9$PM, meaning that selling activity happens and is encouraged by
the highest selling price $\beta^t$ in these slots across the entire
time horizon. Moreover, selling activity results in the peak
conventional generation from $7$PM to $9$PM. Fig.~\ref{fig:costAB}
compares the optimal costs for the two cases. It can be seen that
the optimal costs of conventional generation and worst-case
transaction of Case B are higher than those of Case A, which can be
explained by the higher transaction prices and the resultant larger
DG output for Case B.

The optimal power scheduling of class-2 elastic load is depicted in Fig.~\ref{fig:PEqA} for Case A.
Due to the start time $S_q$ (cf.~Table~\ref{tab:loadK2}), zero power is scheduled for the
class-2 load 1, 3, and 4 from $4$PM to $6$PM while from $4$PM to $7$PM for the load 2.
The decreasing trend for all such loads is due to the decreasing weights $\{\pi_q^t\}$
from $S_q$ to $T_q$, which is established from the fast charging
motivation for the PHEVs, for example.

Figs.~\ref{fig:PBjB} depicts the optimal charging or discharging power of the DSs for Case B.
Clearly, all DSs are discharging during the three slots of $7$PM, $8$PM, and $9$PM. This results
from the motivation of selling more or purchasing less power because both purchase
and selling prices are very high during these slots (cf.~Table~\ref{tab:price}).
The charging (discharging) activity can also be reflected by the stored energy of the battery devices
shown in Fig.~\ref{fig:BjB}. Note that, starting from the initial energy $5$kWh at $4$PM,
the optimal stored energy of all units are scheduled to have $5$kWh at $12$AM, which
satisfies the minimum stored energy requirement for the next round of scheduling time horizons.

%The coincidence of stored energy and the battery capacity for the several last slots is the result of the postulated battery
%cost penalizing the deviation from capacity. The increasing stored energy
%behavior can be explained by the increasing weights $\{\varpi_j^t\}$.
%and the persistent charging behavior (positive $P_{B_j}^t$).

Finally, Fig.~\ref{fig:ratiocostB} shows the effect of different selling
prices $\{\beta^t\}$ on the optimal energy costs, where Case B is
studied with fixed purchase prices $\{\alpha^t\}$. It can be clearly
seen that the net cost decreases with the increase of the
selling-to-purchase-price ratio $\beta^t/ \alpha^t$. When this ratio increases,
the microgrid has a higher margin for revenue from the transaction mechanism,
which yields the reduced worst-case transaction cost.
%Finally, it is worth mentioning that if a similar
%numerical test was implemented for Case A, but no significant change
%emerged in either type of cost. The reason is that the low real-time
%prices $\{\alpha^t\}$ and $\{\beta^t\}$ in Case A do not
%considerably affect the transaction amount for different selling
%prices .

\section{Conclusions and Future Work}\label{sec:Conclusions}
A distributed energy management approach was developed tailored for
microgrids with high penetration of renewable energy sources. By
introducing the notion of committed renewable energy, a novel model
was introduced to deal with the challenging constraint of the
supply-demand balance raised by the intermittent nature of renewable
energy sources. Not only the conventional generation costs,
utilities of the adjustable loads, and distributed storage costs
were accounted for, but also the worst-case transaction cost was
included in the objective. To schedule power in a distributed
fashion, the dual decomposition method was utilized to decompose the
original problem into smaller subproblems solved by the LCs of
conventional generators, dispatchable loads, DS units and the RES.
%The numerical tests demonstrated the optimal power generated,
%consumed by and delivered to the DERs, and the controllable loads
%across the entire time horizon.

A number of interesting research directions open up towards
extending the model and approach proposed in this paper. Some
classical but fundamental problems, such as the optimal power flow
(OPF) and the unit commitment (UC) problems are worth
re-investigating with the envisaged growth of RES usage in microgrids.

%For example, OPF is a non-convex hard problem, whose
%globally optimal solution can be attained for radial distribution
%networks~\cite{LavaeiL11}. Effort must also be made towards
%distributed OPF with the high penetration of RES. Likewise, UC
%problems which are essentially mixed integer programs, are worthy of
%further study.

\appendices

\section{Enhancing the Bundle Method}\label{Appe:Bundle}
Using an auxiliary variable $r$, \eqref{prob:bundle} can be re-written as
\begin{subequations}
\label{prob:QPprimal}
\begin{align}
&\min\limits_{\mathbf{p}, r}~\quad r+\frac{\rho_{\ell}}{2}\|\mathbf{p}-\mathbf{y}_{\ell}\|^2 \\
&\st~\quad \tilde{G}(\mathbf{p}_{i})+\mathbf{g}_{i}^{\prime}(\mathbf{p}-\mathbf{p}_{i}) \le r,~i = 0,1,\ldots,\ell.
\end{align}
\end{subequations}
Introducing multipliers $\bm{\xi} \in \mathbb{R}_{+}^{\ell+1}$, the Lagrangian is given as
\begin{align}
\label{QPLagrang}
\mathcal{L}(r, \mathbf{p}, \bm{\xi}) = \bigg(1-\sum_{i=0}^{\ell+1}\xi_i\bigg)r+\frac{\rho_{\ell}}{2}\|\mathbf{p}-\mathbf{y}_{\ell}\|^2 \nonumber \\
+\sum_{i=0}^{\ell+1}\xi_i\big(\tilde{G}(\mathbf{p}_{i})+\mathbf{g}_{i}^{\prime}(\mathbf{p}-\mathbf{p}_{i})\big).
\end{align}
Optimality condition on $\mathbf{p}$, i.e.,
$\nabla_\mathbf{p}\mathcal{L}(r, \mathbf{p}, \bm{\xi})=\mathbf{0}$ ,
yields
\begin{align}
\label{OptP}
\mathbf{p}^{*} =  \mathbf{y}_{\ell}-\frac{1}{\rho_{\ell}}\sum_{i=0}^{\ell+1}\xi_i\mathbf{g}_{i}.
\end{align}
Substituting \eqref{OptP} into \eqref{QPLagrang}, the dual of
\eqref{prob:QPprimal} is
\begin{subequations}
\label{QPdual}
\begin{align}
&\hspace{-4.0mm} \max\limits_{\bm{\xi}} \,  -\frac{1}{2\rho_{\ell}}\left\|\sum_{i=0}^{\ell+1}\xi_i\mathbf{g}_{i}\right\|^2 +
\sum_{i=0}^{\ell+1}\xi_i\big(\tilde{G}(\mathbf{p}_{i})+\mathbf{g}_{i}^{\prime}(\mathbf{y}_{\ell}-\mathbf{p}_{i})\big)\\
&\st ~\quad  \bm{\xi}\succeq \mathbf{0},~\mathbf{1}^{\prime}\bm{\xi}= 1
\end{align}
\end{subequations}
where $\mathbf{1}$ is the all-ones vector.

Note that~\eqref{QPdual} is essentially a QP over the simplex in
$\mathbb{R}^{\ell+1}$, which can be solved very efficiently.
%It is also worth noting that some enhanced variants of the bundle method, such as dynamic update
%of proximity weight $\rho_{\ell}$, linear search, and Moreau-Yosida regularization may also be
%incorporated~\cite{Kiwiel90,Lemarechal97}.

\section{Proofs of Propositions}
To prove Propositions~\ref{prop:Vchara} and~\ref{prop:Vcomb}, the
following lemma is needed, which shows sufficient and necessary
conditions for a point to be a vertex of a polytope represented as a
linear system~\cite[Sec.~3.5]{Faigle}.
\begin{lemma}
\label{lemma:Vchara} For a polytope $\cP :=\{\mathbf{x} \in
\mathbb{R}^n|\mathbf{A}\mathbf{x} \preceq \mathbf{c}\}$, a point
$\mathbf{v} \in \cP$ is a vertex if and only if there exists a
subsystem $\tilde{\mathbf{A}}\mathbf{x} \preceq \tilde{\mathbf{c}}$
of $\mathbf{A}\mathbf{x} \preceq \mathbf{c}$ so that
$\textrm{rank}(\tilde{\mathbf{A}})=n$ and $\mathbf{v}$ is the unique
(feasible) solution of $\tilde{\mathbf{A}}\mathbf{v} =
\tilde{\mathbf{c}}$.
\end{lemma}

\subsection{Proof of Proposition~\ref{prop:Vchara}}\label{Appe:Vchara}
The polytope $\cA := \{\mathbf{a} \in \mathbb{R}^n|\underline{\mathbf{a}}\preceq \mathbf{a} \preceq
\overline{\mathbf{a}}, a^{\min}\le  \mathbf{1}^{\prime}\mathbf{a}
\le a^{\max}\}$ can be re-written as $\cA := \{\mathbf{a} \in
\mathbb{R}^n|\mathbf{A}\mathbf{a}\preceq \mathbf{c}\}$, where
$\mathbf{A}:=[\mathbf{I}_{n \times n}, -\mathbf{I}_{n \times n},
\mathbf{1},-\mathbf{1}]^{\prime}$ and
$\mathbf{c}:=[\overline{\mathbf{a}}^{\prime},-\underline{\mathbf{a}}^{\prime},a^{\max},-a^{\min}]^{\prime}$.
By Lemma~\ref{lemma:Vchara}, enumerating vertices of $\cA$ is
equivalent to finding all feasible solutions of the linear
subsystems $\tilde{\mathbf{A}}\mathbf{a} = \tilde{\mathbf{c}}$, such
that rank-$n$ matrix $\tilde{\mathbf{A}}$ is constructed by
extracting rows of $\mathbf{A}$. It can be seen that such full
column-rank matrix $\tilde{\mathbf{A}}$ can only have two forms
(with row permutation if necessary): i) $\tilde{\mathbf{A}}_1 =
\textrm{diag}(\mathbf{d})$ with $d_i \in \{-1,1\},~i=1,\ldots,n$;
ii) $\tilde{\mathbf{A}}_2(i,:) =
\pm\mathbf{1}^{\prime},~i=1,\ldots,n$, and
$\tilde{\mathbf{A}}_2(j,:) = \tilde{\mathbf{A}}_1(j,:),~\forall
j\neq i$. Basically, $\tilde{\mathbf{A}}_1$ is constructed by
choosing $n$ vectors as a basis of $\mathbb{R}^n$ from the first
$2n$ rows of $\mathbf{A}$. Substituting any row of
$\tilde{\mathbf{A}}_1$ with $\pm\mathbf{1}^{\prime}$, forms
$\tilde{\mathbf{A}}_2$. Finally, by solving all the linear
subsystems of the form $\tilde{\mathbf{A}}_{k}\mathbf{a} =
\tilde{\mathbf{c}}_{k}$, for $k=1,2$, Proposition~\ref{prop:Vchara}
follows readily.

\subsection{Proof of Proposition~\ref{prop:Vcomb}}\label{Appe:Vcomb}
The polytope $\cB := \{\mathbf{b}\in \mathbb{R}^n|\underline{\mathbf{b}}\preceq \mathbf{b}\preceq
\overline{\mathbf{b}}, b^{\min}_s \le \mathbf{1}_{n_s}^{\prime}\mathbf{b}_s \le b^{\max}_s, s =
1,\ldots,S\}$ can be re-written as $\cB := \{\mathbf{b} \in
\mathbb{R}^n|\mathbf{B}\mathbf{b}\preceq \mathbf{c}\}$, where
$\mathbf{B} := \textrm{diag}(\mathbf{B}_1,\ldots,\mathbf{B}_S)$,
$\mathbf{c} :=
[\mathbf{c}_1^{\prime},\ldots,\mathbf{c}_S^{\prime}]^{\prime}$,
$\mathbf{B}_s:=[\mathbf{I}_{n_s \times n_s}, -\mathbf{I}_{n_s \times
n_s}, \mathbf{1},-\mathbf{1}]^{\prime}$, and
$\mathbf{c}_s:=[\overline{\mathbf{a}}_s^{\prime},-\underline{\mathbf{a}}_s^{\prime},b_s^{\max},-b_s^{\min}]^{\prime}$
for $s=1,\ldots,S$.

Similarly by Lemma~\ref{lemma:Vchara}, all the vertices of $\cB$ can
be enumerated by solving $\tilde{\mathbf{B}}\mathbf{b} =
\tilde{\mathbf{c}}$, where the rank-$n$ matrix $\tilde{\mathbf{B}}$
is formed by extracting rows of $\mathbf{B}$. Due to the block
diagonal structure of $\mathbf{B}$, it can be seen that the only way
to find its $n$ linear independent rows is to find $n_s$ linear
independent vectors from the rows corresponding to $\mathbf{B}_s$
for $s=1,\ldots,S$. In other words, the vertices $\mathbf{b}^{\tV}$
can be obtained by concatenating all the individual vertices
$\mathbf{b}_s$ as stated in Proposition~\ref{prop:Vcomb}.

%\IEEEtriggercmd{\enlargethispage{-5.35in}}
\IEEEtriggeratref{30}

%%%%%%%%%%%%%%%%%%%%%%%%%%%%%%%%%%%%%%%%%%%%%%
\bibliographystyle{IEEEtran}
\bibliography{biblio}

\end{document}